\documentclass[11pt, a4paper]{article}

\usepackage[activate={true,nocompatibility},final,tracking=true,kerning=true,spacing=true,factor=1100,stretch=10,shrink=10]{microtype} 

\SetProtrusion{encoding={*},family={bch},series={*},size={6,7}}
              {1={ ,750},2={ ,500},3={ ,500},4={ ,500},5={ ,500},
               6={ ,500},7={ ,600},8={ ,500},9={ ,500},0={ ,500}}
               
\SetExtraKerning[unit=space]
    {encoding={*}, family={bch}, series={*}, size={footnotesize,small,normalsize}}
    {\textendash={400,400}, 
     "28={ ,150}, 
     "29={150, }, 
     \textquotedblleft={ ,150}, 
     \textquotedblright={150, }} 
     
\SetExtraKerning[unit=space]
   {encoding={*}, family={qhv}, series={b}, size={large,Large}}
   {1={-200,-200}, 
    \textendash={400,400}}

\usepackage[usenames,dvipsnames]{xcolor} 

\usepackage[left=3.5cm, right=3.5cm, top=3.5cm, bottom=3.5cm, showframe=false]{geometry}
\usepackage[affil-it]{authblk}
\usepackage[fleqn]{amsmath}
\usepackage{amsthm}
\usepackage{amssymb}
\usepackage{booktabs} 
\usepackage{graphicx}
\usepackage{bbm} 
\usepackage{here}
\usepackage{subfig}
\usepackage[round]{natbib}
\usepackage{hyperref}
\usepackage{titlesec}
\usepackage{thmtools}
 \usepackage{setspace}
\usepackage{lmodern}
	\usepackage{dsfont}
\usepackage{lscape}
 \usepackage[latin1]{inputenc} 
 \usepackage[english]{babel}
\usepackage{charter} 
\usepackage[T1]{fontenc} 

					\newcommand{\annmathinform}{{Annales Mathematicae et Informaticae}} 
					
					\newcommand{\arxiv}{{arXiv}}
					
					\newcommand{\computationalstatisticsanddataanalysis}{{Computational Statistics and Data Analysis}}
			
					\newcommand{\isrnjournal}{{ISRN Probability and Statistics}}

					\newcommand{\jasa}{{Journal of the American Statistical Association}} 
					
					\newcommand{\journalofeconometrics}{{Journal of Econometrics}} 
					
					\newcommand{\journalmultvaranal}{{Journal of Multivariate Analysis}} 
					
					\newcommand{\journalroyalsociety}{{Journal of the Royal Statistical Society: Series B}} 
					
					\newcommand{\journaltimeseranal}{{Journal of Time Series Analysis}}

					\newcommand{\metrika}{{Metrika}}
					\newcommand{\probtheoryrelatedfields}{{Probability Theory and Related Fields}}

					\newcommand{\scandinavian}{{Scandinavian Journal of Statistics}}
					
					\newcommand{\statisticasinica}{{Statistica Sinica}}
					
					\newcommand{\stochprocappl}{{Stochastic Processes and their Applications}}
					
					\newcommand{\statprob}{{Statistics and Probability Letters}}

					\newcommand{\theorprobabappl}{{Theory of Probability and Its Applications}}
					
					\newcommand{\annmathstat}{{The Annals of Mathematical Statistics}}
					 
					\newcommand{\annstat}{{The Annals of Statistics}}
					
					\newcommand{\annprob}{{The Annals of Probability}}
					
					\newcommand{\journaltest}{{TEST}}
					
\titleformat*{\section}{\Large\bfseries}
\titleformat*{\subsection}{\large\bfseries}
\titleformat*{\subsubsection}{\bfseries}
\titleformat*{\paragraph}{\bfseries}

\newcommand{\cA}{{\cal A}}
\newcommand{\cB}{{\cal B}}
\newcommand{\cC}{{\cal C}}

\newcommand{\cO}{{\mathcal O}}

\newcommand{\mN}{\mathbb{N}}
\newcommand{\mR}{\mathbb{R}}
\newcommand{\mZ}{\mathbb{Z}}

\newcommand{\bB}{\mathbf{B}}

\newcommand{\bW}{\mathbf{W}}

\newcommand{\bX}{\mathbf{X}}

\newcommand{\bY}{\mathbf{Y}}

\newcommand{\lpm}{$L^\kappa$-$m$}
\newcommand{\lfourm}{$L^4$-$m$}

\newcommand{\bSigma}{{\mathbf{\Sigma}}}

\newcommand{\bmu}{\boldsymbol{\mu}}
\newcommand{\bDelta}{\boldsymbol{\Delta}}

\newcommand{\Cov}{{\rm Cov}}
\newcommand{\diag}{{\rm diag}}
\newcommand{\Var}{{\rm Var}}

\newcommand{\ntheta}{{\lfloor n \theta \rfloor}}

\theoremstyle{plain}
\newtheorem{theorem}{Theorem}[section]
\newtheorem{proposition}[theorem]{Proposition}
\newtheorem{corollary}[theorem]{Corollary} 

\theoremstyle{definition}
\newtheorem{definition}[theorem]{Definition}
\newtheorem{notations}{Notation}  

\theoremstyle{remark}
\newtheorem{remark}[theorem]{Remark}

{\bf}{\it}
{\bf}{\it}
{\bf}{\it}
{\bf}{\it}
{\bf}{\it}
{\bf}{\it}
{\bf}{\it}
{\bf}{\it}
{\bf}{\it}
{\bf}{\it}

\declaretheoremstyle[
    spaceabove=0.5cm, 
    spacebelow=0.5cm, 
    headfont=\normalfont\bfseries,
    notefont=\normalfont\bfseries, 
    notebraces={}{}, 
    bodyfont=\normalfont, 
    headpunct={.} 
]{myassumption}

	  	\theoremstyle{myassumption}
		\newtheorem{taggedassumptionX}{Assumption} 
		\newenvironment{taggedassumption}[1]
		 {\renewcommand\thetaggedassumptionX{#1}\taggedassumptionX}
 		{\endtaggedassumptionX}

\makeatletter
\def\@maketitle{%
  \newpage
  \null
  \vskip 2em%
  \begin{center}%
  \let \footnote \thanks
    {\Large\bfseries \@title \par}%
    \vskip 1.5em%
    {\normalsize
      \lineskip .5em%
      \begin{tabular}[t]{c}%
        \@author
      \end{tabular}\par}%
    \vskip 1em%
    {\normalsize \@date}%
  \end{center}%
  \par
  \vskip 1.5em}
\makeatother

		\colorlet{linkcolour}{cyan!50!blue}
		\colorlet{urlcolour}{magenta}
		\hypersetup{colorlinks=true, linkcolor=linkcolour, citecolor=linkcolour, urlcolor=urlcolour}
		
\numberwithin{equation}{section}

\begin{document}

\title{A Darling-Erd{\H o}s-type CUSUM-procedure \\ for functional data II}

\author{Leonid Torgovitski\thanks{E-mail: \texttt{ltorgovi@math.uni-koeln.de}} \thanks{Research partially supported by the Friedrich Ebert Foundation.}}
\affil{Mathematical Institute, University of Cologne\\ Weyertal 86-90, 50931, Cologne, Germany}
\date{}

\maketitle
{\footnotesize
\begin{abstract}
{\footnotesize
This article considers testing for mean-level shifts in functional data. The class of the famous Darling-Erd\H{o}s-type cumulative sums (CUSUM) procedures is extended to functional time series under short range dependence conditions which are satisfied by functional analogues of many popular time series models including the linear functional AR and the non-linear functional ARCH. We follow a data driven, projection-based approach where the lower-dimensional subspace is determined by (long run) functional principal components which are eigenfunctions of the long run covariance operator. This second-order structure is generally unknown and estimation is crucial - it plays an even more important role than in the classical univariate setup because it generates the finite-dimensional subspaces. We discuss suitable estimates and demonstrate empirically that altogether this change-point procedure performs well under moderate temporal dependence. 

Moreover, Darling-Erd\H{o}s-type change-point estimates based on (long run) functional principal components as well as the corresponding {\frqq fully-functional\flqq} counterparts are provided and the testing procedure is finally applied to publicly accessible electricity data from a German power company. 
}
\end{abstract}

 }
\vspace{1.5cm}
\section{Introduction}

The interest and the research activities in \textit{\frqq change-point analysis\flqq} for multivariate, high-dimensional and especially for functional data are enormous which is a consequence of the increasingly growing computational capacities. These activities are reflected by the amount and the high frequency of published works and in particular by survey articles that appeared recently.\footnote{Cf., e.g., \citet{kokoszka2012}, \citet{aue2013} and the invited discussion paper by \citet{test2014}.} One of the fundamental and most studied problems in change-point analysis is concerned with a simple abrupt change in the mean - the \textit{\frqq at most one change\flqq} (AMOC) model.  
\begin{itemize}
\item We consider this problem in the functional setup, i.e.\ where each observation is a curve and the mean is a deterministic function. 
\item We want to know whether the overall shapes of the mean-curves have changed over time at some arbitrary time point or not. 
\end{itemize}
Our investigations are based on the work of \citet{berkes2009functional} who studied the same problem and introduced a (differently weighted) nonparametric CUSUM procedure for detection of changes in the mean of functional observations in the i.i.d.\ setting. \citet{berkes2009functional} suggested an intuitive approach which essentially relies on a multivariate CUSUM by projecting the functional time series on a finite dimensional subspace which captures the dynamics of the data in a beneficial manner in order to obtain reasonable power. The authors proposed to select the subspace spanned by functional principal components (FPC's), i.e.\ the eigenfunctions of the covariance operator. This approach is motivated by their well known optimality properties regarding dimension reduction (cf. \citet{ramsay2005functional}). Since then, FPC's - which play widely known an outstanding role in functional data analysis - have been successfully incorporated into many further functional \textit{\frqq stability-testing\flqq} procedures under independence as well as more recently under dependence (cf. \citet{horvath2012fda}  for an overview and also \citet{berkes2009functional}, \citet{hoermann2010weakly} and \citet{kirch2012functional} in particular for the change in the mean setting). Later on, it was realized that long run FPC's, given as eigenfunctions of the so-called long run covariance operator, are advantageous (cf. \citet{reeder2012twoB, horvath2013testing} and \citet{torg2014}).

In this article we will stick to the latter approach. We pick up and continue the work of \citet{torg2014} (cf. also \citet{zhou2011maximum}) extending the results from the $m$-dependent setting to the more challenging and realistic models of weakly dependent time series with a focus on the framework of \lpm-approximability (cf., e.g.,  \citet{hoermann2010weakly} and \citet{reeder2012twoB}). As in \citet{torg2014}, the procedure will be based on the dimension-reduction approach of \citet{berkes2009functional}. To establish asymptotics we will incorporate several steps outlined in \citet{torg2014} and combine them with results of \citet{berkes2011split, berkes2013rice}, \citet{reeder2012twoB} and \citet{Aue2014}.

This article contributes to the massive amount of recent works on change-point testing and estimation in functional (or generally high-dimensional) data and is a revised version of \citet{torg2014approx}. 
\begin{enumerate}
\item On one hand, our results on long run covariance operator estimation {\frqq complement\flqq} the findings of \citet{1lukasz2014}.\footnote{Cf., also \citet{whipple2014} and \citet{berkes2015rice}.} Their results are slightly stronger but our proof technique is different and thus is of separate interest.
\item On the other hand, our results on the (multivariate) Darling-Erd{\H o}s limit theorems complement the related discussion of \citet{Kirch2014}. Here, we show additional conditions that emerge due to dimension reduction, i.e.\ due to the transition from the functional to the multivariate settings. Moreover, we verify conditions for the multivariate Darling-Erd{\H o}s asymptotics explicitly under the specific dependence concept of \lpm-approximability.  
\item Furthermore, we discuss the relation of projection-based and fully-functional estimates for change-points. This contributes to the investigation of some related estimates in the recent work of \citet{aue2015}. 
\item Finally, we demonstrate the performance of the Darling-Erd{\H o}s-type CUSUM procedure in Monte Carlo simulations and conduct an analysis of a real-life {\frqq electricity dataset\flqq}. Note that the {\frqq synthetic\flqq} simulations presented here and in the previous version \citet{torg2014approx} are used for comparison by \citet{tewes2015}.
\end{enumerate}
 
\begin{notations}
In order to formalize the testing problem we have to introduce some notation first.  We consider \textit{functional data} $X(\cdot)$ as a \textit{random element} on some probability space $(\Omega, \cA, P)$ with the state space $L^2[0,1]$. 
Throughout, let $L^2[0,1]$ denote the space of square-integrable functions with respect to the Lebesgue measure on $[0,1]$ equipped with the usual inner product and the corresponding norm, denoted by $\langle v,w\rangle=\int\!v(t)w(t)dt$ or $\|v\|$ for $v,w\in L^2[0,1]$, respectively. We also assume product measurability of $X(t)=X(t,\omega)$ with respect to $(t,\omega)\in[0,1]\times \Omega$. The mean of $X(\cdot)$ is defined as the unique function $\mu(\cdot)$, such that $\int\!x(t)\mu(t)dt=E\int\!x(t)X(t)dt$ holds true for all $x\in L^2[0,1]$ given that $E\|X\|<\infty$. 
\end{notations} 

We assume that the observable sequence $\{X_i(\cdot)\}_{i\in\mZ}$ consists of $L^2[0,1]$-valued random elements which are given by the functional \textit{\frqq signal plus noise\flqq} model
\begin{equation}\label{eq:decomp_innov_F}
     X_i(t) = \mu_i(t) + Y_i(t),\qquad t\in[0,1],
\end{equation}
with mean functions $\mu_i(\cdot)\in L^2[0,1]$ and with innovations fulfilling our basic \autoref{ass:assumptionM}, below. The dependence structure of the innovations will be specified later on.  We want to test retrospectively the null hypothesis of no change in the mean, i.e.
\begin{align*}
	H_0: \qquad &\mu_1(\cdot)=\ldots = \mu_n(\cdot)\\
\intertext{against the alternative of a change in the mean}
 	H_A:  \qquad &\mu_1(\cdot)= \ldots =  \mu_{\lfloor n \theta \rfloor}(\cdot) \neq  \mu_{\lfloor n \theta \rfloor+1}(\cdot)=\ldots= \mu_n(\cdot)
\end{align*}
at some unknown time point characterized by some (unknown) constant change parameter $\theta\in (0,1)$.

\begin{taggedassumption}{(M)}\label{ass:assumptionM} \hfill
\begin{itemize}
\item[(i)] The functional innovation sequence $\{Y_i\}_{i\in \mZ}$ is centered and strictly stationary; 
\item[(ii)] $E\|Y_1\|^\nu<\infty$ for some $\nu>2$.
\end{itemize}
\end{taggedassumption}

\textit{The structure of this article is as follows:} In \autoref{sect:weakdep} we formulate the dependence concept of \lpm-approximability for our observations. In \autoref{sect:p2test_procedure} we present the testing procedure together with the asymptotics under the null hypothesis and under the alternative. \autoref{sect:estimation} focuses on estimation of long run FPC's. The performance is finally demonstrated in \autoref{sect:simulAX} including an application example. All proofs are postponed to \autoref{sect:proofII}.
 
\section{Weakly dependent time series}\label{sect:weakdep}

We consider the \textit{\frqq mathematically convenient\flqq} concept of \lpm-approximable time series which is currently of major interest in univariate, multivariate and functional settings and covers many relevant time series models (cf., e.g., amongst many others \citet{Aue2009}, \citet{hoermann2010weakly}, \citet{kirch2012functional}, \citet{reeder2012twoB, horvath2013testing}, \citet{jirak2012, jirak2013}, \citet{berkes2013rice}, \citet{chochola2013} and \citet{1lukasz2014}).  \citet{hoermann2010weakly} and \citet{berkes2011split} contain extensive discussions and comparisons to other related dependence concepts.

We formulate the dependence concept in general real and separable Hilbert spaces $H$ but having the special cases $H=L^2[0,1]$ and $H=\mR^d$ in mind. The reason for this is that we will use the notion of \lpm-approximability in both spaces because we will also deal with appropriate $\mR^d$-valued approximations of the original $L^2[0,1]$ valued time series in our proofs. For a moment, let $\|\cdot\|_H$ denote a norm in the space $H$ and recall that $E(\|X\|_H^\kappa)^{1/\kappa}$ is the $L^\kappa(\Omega, P)$-norm for the real-valued random variable $\|X\|_H$. Later on we will write $\|\cdot\|$ for the $L^2$ norm or $|\cdot|$ for the Euclidean norm, respectively.
\begin{definition}\label{def:lpm}
The $H$-valued sequence $\{Y_i\}_{i\in \mZ}$, defined on some common probability space $(\Omega, \cA, P)$,  is \lpm-approximable with rate $\delta(m)=o(1)$, $\delta(m)\geq 0$  and with $\kappa\geq 2$ iff $E[\|Y_0\|_H^{\kappa}]<\infty$ and the following conditions hold true:
\begin{itemize}
\item[1.] The $Y_i$'s admit a \textit{\frqq Bernoulli shift\flqq} representation, i.e.
\begin{equation}\label{eq:non_lin_modelp2}
	Y_i=f(\ldots, \varepsilon_{i+1}, \varepsilon_i, \varepsilon_{i-1}, \ldots),
\end{equation}
where the innovations $\varepsilon_i$ are i.i.d.\ $S$-valued random elements, $S$ is some measurable space and $f$ is a measurable mapping $f: S^\infty \rightarrow H$.
\item[2.]  The $Y_i$'s are approximated by $m$-dependent random variables in the sense that, as $m\rightarrow \infty$,  
\begin{equation}\label{eq:F_ass_decay}
	E(\|Y_0-Y_0^{(m)}\|^\kappa_H)^{1/\kappa}= \cO(\delta(m)),
\end{equation}
where $Y_i^{(m)}$ are $m$-dependent copies of $Y_i$ defined via
\begin{equation}\label{eq:F_mcopy}
Y^{(m)}_i= f(\ldots, \varepsilon^{(m,i)}_{i+M}, \varepsilon_{i+(M-1)}, \ldots, \varepsilon_{i} ,\ldots , \varepsilon_{i-(M-1)}, \varepsilon^{(m,i)}_{i-M},\ldots)
\end{equation}
with $M=\max\{\lfloor m/2 \rfloor,1\}$ for all integer $m\geq 0$ and where the family 
\[
	\{\varepsilon_r, \varepsilon_{i}^{(k,j)}, \; i,j,r,k\in \mZ,\; k\geq 0\}
\]
is i.i.d.\
\item[3.] If $\{Y_i\}$ is function-valued with $H=L^2[0,1]$ then $(\cB_{[0,1]} \otimes \cA)-\cB_{\mR}$ measurability is assumed.
\end{itemize}
\end{definition}Typical conditions on the rate $\delta(m)$ are summability $\sum_{m=1}^\infty\delta(m)<\infty$, polynomial decay $\delta(m)=\cO(m^{-\nu})$ for some $\nu>2$ or exponential decay $\delta(m)=\cO(\exp(-cm))$ with some $c>0$, where the latter is here the strongest condition but already satisfied for many models such as, e.g., the $H$-valued AR(1). 

\begin{remark}We are interested in the causal case of \lpm-approximability, i.e.\ that $Y_i=f(\varepsilon_i, \varepsilon_{i-1}, \ldots)$ holds true. This is a special case of \eqref{eq:non_lin_modelp2} but the two-sided {\frqq noncausal\flqq} representation \eqref{eq:non_lin_modelp2} appears to be useful in the proofs, where we will deal with time-inversed \lpm-approximable time series $\{Y_{-i}\}_{i\in \mZ}$. Therefore, observe that 
\[
	Y_{-i}=f(\varepsilon_{-i}, \varepsilon_{-(i+1)}, \varepsilon_{-(i+2)}, \ldots)
\]
is noncausal but still \lpm-approximable according to the two-sided Definition \ref{def:lpm}, above.
\end{remark}
\begin{remark}\label{rem:different_cond}
Note that some recent literature (e.g.\ \citet{berkes2013rice} and \citet{horvath2013testing}) works with a slightly modified condition \eqref{eq:F_ass_decay} where the left-hand side of \eqref{eq:F_ass_decay} is substituted by $\smash{E(\|Y_0-Y_0^{(m)}\|^{\kappa}_H)^{1/\tilde{\kappa}}}$ with some $\tilde{\kappa}>\kappa$.
\end{remark}

\section{The testing procedure}\label{sect:p2test_procedure}

For the sake of generality and clarity, the testing procedure will be described in a unifying functional framework where we separate and highlight those conditions which essentially allow us to derive suitable asymptotics without a priori specifying a particular time series model or a dependence concept. The conditions presented below in this section were (to some degree) implicitly contained in \citet{torg2014}. Here, the theoretical focus will be more on verification of all stated conditions for \lpm-approximable time series.

The CUSUM procedure, which will be presented below, belongs to the class of \textit{\frqq FPC-based approaches\flqq} and utilizes the second-order structure of the time series for an appropriate data-driven subspace selection. As mentioned in the introduction, we will assume a functional time series $\{X_i\}$ with functional innovations $\{Y_i\}$ and work with long run FPC's following, e.g.,  \citet{reeder2012twoB, horvath2013testing} and \citet{torg2014}. Those are eigenfunctions of the long run covariance operator $\cC$ of $\{Y_i\}_{i\in \mZ}$ which, given \autoref{ass:assumptionM}, can be formally defined as an integral operator
\[
	\big(\cC x\big)(t)= \int{\!\zeta(t,s) x(s)} ds,
\]
$t\in[0,1]$, $x \in L^2[0,1]$, $\int:=\int_0^1$ with kernel
\begin{equation}\label{eq:kernel_zeta_F}
	\zeta(t,s)=\sum_{r\in\mZ} E[Y_0(t) Y_r(s)].
\end{equation}
This operator is well-defined if
\begin{equation}\label{eq:kernelL2}
	\zeta \in L^2([0,1]\times[0,1])
\end{equation}
holds true in which case it is symmetric Hilbert-Schmidt and positive. Hence, $\cC$ can be written using the spectral decomposition as 
\begin{equation}\label{eq:cov_series}
	(\cC x)(t) = \sum_{j=1}^{\infty}\lambda_{j}\Big[\int\! x(s)v_{j}(s)ds\Big] v_{j}(t), 
\end{equation}
$t\in[0,1]$, $x \in L^2[0,1]$, with real eigenvalues $\lambda_{1}\geq \lambda_{2} \geq \ldots \geq 0$ in descending order and corresponding orthonormal eigenfunctions $v_1(t),v_2(t),\ldots$. The convergence of the series \eqref{eq:kernel_zeta_F} and \eqref{eq:cov_series} above is meant in the $L^2([0,1]\times[0,1])$ or $L^2[0,1]$ sense, respectively. Clearly, the eigenvalues are non-negative due to the positiveness of $\cC$.

\begin{remark} Notice, that the decomposition \eqref{eq:cov_series} is obviously ambiguous. One reason are signs: each eigenvalue $\lambda_j$ is always associated at least with two eigenfunctions $v_{j}$ and $-v_{j}$; Another source of ambiguity is the geometric multiplicity of $\lambda_j$ being larger than one. The statistic will be constructed in such a way that it will turn out to be invariant under changes of signs. Also, conditions will be imposed (cf. \eqref{eq:eigenvalues_distinct} below)  to ensure that relevant eigenspaces are only one-dimensional - both is standard practice for FPC-based statistics. (In the following we tacitly restrict ourselves to $\lambda_d>0$, i.e.\ the functional data is at least $d$-dimensional.) 
\end{remark} 

For testing of $H_0$ against $H_A$ we will work with the following CUSUM statistic 
\begin{align}
\begin{split}\label{eq:ml_stat_F}
	T_n(X;v,\lambda)=\max_{1 \leq  k < n}w(k/n) \left(\frac{\eta^2_{k,1}}{\lambda_1}+\ldots +\frac{\eta^2_{k,d}}{\lambda_d}\right)^{1/2}
\end{split}
\end{align}
with scores
\[
	\eta_{k,r}=\eta_{k,r}(X;v)=n^{-1/2} \int\!\sum_{i=1}^{k}[X_i(t) - \bar{X}_{n}(t)]v_r(t)dt
\]
where $w(t)=(t(1-t))^{-1/2}$ is a suitable weight function related to the variance of a Brownian bridge and $d\in\mN$ is a fixed positive integer specifying the dimension of the subspace chosen for projecting.  Notice, that in \eqref{eq:ml_stat_F} $X$, $v$ or $\lambda$ represent $\{X_i\}_{i\in \mN}$, $\{v_i\}_{i\in \mN}$ or $\{\lambda_i\}_{i\in \mN}$, respectively. The right-hand side of \eqref{eq:ml_stat_F} can also be written compactly using vector-notation as
\[
	T_n(X;v,\lambda)=\max_{1 \leq  k < n} w(k/n) |n^{-1/2}\sum_{i=1}^{k}(\bX_i - \bar{\bX}_n)|_\bSigma,
\]
with the norm $|\cdot|_{\bSigma}=|\bSigma^{-1/2}\cdot|$ and where $\bSigma(\lambda)=\diag(\lambda_1,\ldots, \lambda_d)$  is a proper standardization matrix (cf. \citet{torg2014} and \autoref{rem:covmatrix} below). This presentation emphasizes that our CUSUM is based on the multivariate {\frqq projected version\flqq} of the data \eqref{eq:decomp_innov_F}, i.e.\ on
\begin{equation}\label{eq:proj_innov}
	\bX_i = \bmu_i + \bY_i,
\end{equation}
where the $r$-th components of these vectors are 
\begin{itemize}
\item the data \textit{scores} $\bX_{i,r}=\int\! X_i(t) v_r(t)dt$, 
\item the innovation \textit{scores}  $\bY_{i,r}=\int\! Y_i(t) v_r(t)dt$ 
\item and the projected means $\bmu_{i,r}=\int\! \mu_i(t) v_r(t) dt$ 
\end{itemize}
for $r=1,\ldots,d$, respectively.  
\begin{remark}\label{rem:covmatrix}The matrix $\bSigma(\lambda)$ is the long run covariance matrix of the projected time series $\{\bY_i\}_{i\in\mZ}$ which can be seen utilizing the orthonormality of the eigenfunctions:
\[
	\sum_{k\in\mZ}E[\bY_{0,i}\bY_{k,j}]=\int\! \left(\int\! \zeta(t,s)v_j(s)ds\right)v_i(t)dt=\lambda_j\delta_{i,j}.
\]
The diagonality of $\bSigma$ is hereby one of the main advantages of working with the long run covariance operator $\cC$. Also, the existence of $\bSigma$ is inherited from the existence of the functional counterpart $\cC$.
\end{remark}

In applications, especially in the functional setup, the covariance structure, e.g.\ in our case the covariance operator $\cC$, will be rarely known. Hence, the associated quantities such as the eigenelements $\{v_i\}_{i\in \mN}$ and $\{\lambda_i\}_{i\in \mN}$, are therefore also usually unknown and have to be estimated. Therefore, let $\{\hat{v}_i\}_{i\in \mN}$ be orthonormal functions which together with non-negative scalars $\{\hat{\lambda}_i\}_{i\in \mN}$, $\hat{\lambda}_1\geq \hat{\lambda}_2 \geq \ldots \geq 0$, denote some generic estimates which will be specified later on. Thus, instead of working with $T_n=T_n(X;v,\lambda)$ we will consider 
\begin{align}\label{eq:stat_est_p2}
\begin{split}
	\hat{T}_n&=T_n(X;\hat{v},\hat{\lambda}) =\max_{1 \leq  k < n}w(k/n)\left(\frac{\hat{\eta}^2_{k,1}}{\hat{\lambda}_1}+\ldots +\frac{\hat{\eta}^2_{k,d}}{\hat{\lambda}_d}\right)^{1/2}
\end{split}
\end{align}
with estimated scores
\[
	\hat{\eta}_{k,r}=\eta_{k,r}(X;\hat{v})=n^{-1/2} \int\!\sum_{i=1}^{k}[X_i(t) - \bar{X}_{n}(t)]\hat{v}_r(t)dt
\]
having in mind that formally $\hat{T}_n=\infty$ if $\hat{\lambda}_r=0$ for some $1\leq r \leq d$. The vector notation in the estimated case is given by $\hat{T}_n=\max_{1 \leq  k < n} w(k/n) |n^{-1/2}\sum_{i=1}^{k}(\hat{\bX}_i - \bar{\hat{\bX}}_n)|_{\hat{\bSigma}}$ where each component of the projected time series $\{\hat{\bX}_i\}$ is $\hat{\bX}_{i,r}=\int\! X_i(t) \hat{v}_r(t)dt$ and where $\hat{\bSigma}=\diag(\hat{\lambda}_1,\ldots,\hat{\lambda}_d)$. A natural estimate for the change-point is given analogously to \eqref{eq:ml_stat_F} and \eqref{eq:stat_est_p2} through
\[
	\hat{k}(\hat{v},\hat{\lambda},d)=\operatorname*{arg\,max}_{1 \leq  k < n}w(k/n)\left(\frac{\hat{\eta}^2_{k,1}}{\hat{\lambda}_1}+\ldots +\frac{\hat{\eta}^2_{k,d}}{\hat{\lambda}_d}\right)^{1/2}.
\] 
We will discuss some related estimates in \autoref{rem:ffunc}, below.
\hfill
\\
\hfill
\\
The following conditions (L), (P1), (P2), (A1), (A2) and (B1), (B2) are the main {\frqq building blocks\flqq} which allow us to prove the limiting distribution of  $\hat{T}_n$ under the null hypothesis and consistency under the alternative. Recall that \autoref{ass:assumptionM} is always tacitly assumed. We proceed with conditions under the null hypothesis where the first \autoref{ass:assumptionL} states the availability of a multivariate Darling-Erd\H{o}s-type limit theorem for $T_n$ which will be a cornerstone for our further considerations. (See \citet{berkes2009functional}, \citet{hoermann2010weakly} and \citet{kirch2012functional} for related CUSUM procedures based on multivariate functional central limit theorems.)
\begin{taggedassumption}{(L)}\label{ass:assumptionL} It holds that, as $n\rightarrow \infty$,
\begin{equation}\label{eq:exactlimit}
	\lim_{n\rightarrow\infty} P\left(a(\log n)T_n(Y;v,\lambda)- b_d(\log n)\leq x\right)=\exp(-2\exp(-x))
\end{equation}
for all $x\in\mR$, where $a(t) =(2\log t)^{1/2}$ and $b_d(t) =2 \log t + (d/2) {\log \log t} - \log \Gamma(d/2)$ denote the well known normalizing functions.  
\end{taggedassumption}
 
In the i.i.d.\ setting \autoref{ass:assumptionL} is immediately implied by \citet[Theorem 1.3.1]{horvath1997limit} if (ii) of \autoref{ass:assumptionM} and \eqref{eq:kernelL2} holds true. For strictly stationary $m$-dependent sequences \autoref{ass:assumptionL} follows analogously using strong invariance principles derived in \citet{steinebach1999mdependent} (cf. \citet{torg2014}). Verification of \autoref{ass:assumptionL} under \lpm-approximability will be carried out further below but has now to be based on strong approximations derived recently in \citet{Aue2014}. Further conditions (e.g.\ of mixing-type) which ensure \eqref{eq:exactlimit} are briefly discussed in \citet{Kirch2014}.
\begin{remark}
It is worth noting, that strictly stationary $m$-dependent sequences, as considered in \citet{steinebach1999mdependent} or in \citet{torg2014}, are generally either not \lpm-approximable or that the rate function is unknown.  As pointed out by \citet[Section 3.1]{berkes2011split}, they do not necessarily possess representation \eqref{eq:non_lin_modelp2}. 
\end{remark}
The following conditions on maxima of weighted (backward) partial sums of the innovations together with the subsequent conditions on rates for $\{\hat{v}_i\}_{i\in \mN}$ and $\{\hat{\lambda}_i\}_{i\in \mN}$ will ensure a proper {\frqq interplay\flqq} between the functional data and the multivariate statistic $\hat{T}_n$.
\begin{taggedassumption}{(P1)}\label{ass:assumptionP1} It holds that, as $n\rightarrow \infty$, 
\begin{align*}
	\max_{1 \leq  k < n}k^{-1/2}\| \sum_{i=1}^k Y_i\| &= \cO_P(g(n)),\\
\max_{1 \leq  k < n}k^{-1/2}\| \sum_{i=1}^k Y_{-i}\| &= \cO_P(g(n)),
\end{align*}
where the rate function $g(n)$ will be specified later on.
\end{taggedassumption}

\begin{taggedassumption}{(A1)}\label{ass:assumptionA1}  Under $H_0$ it holds that, as $n\rightarrow \infty$, 
\[
	\max_{i=1,\ldots,d}|\hat{\lambda}_i-\lambda_i|=o_P((\log\log n)^{-1}).
\]
\end{taggedassumption}

\begin{taggedassumption}{(A2)}\label{ass:assumptionA2} Under $H_0$ it holds that, as $n\rightarrow \infty$, 
\[
	\max_{i=1,\ldots,d}|\hat{v}_i - s_i v_i|=o_P((\log\log n)^{-1/2}/g(n)),
\]
where $s_i$'s are random with $s_i\in\{1,-1\}$ and the rate function $g(n)$ is the same as in \autoref{ass:assumptionP1}.
\end{taggedassumption}
The random $s_i$'s are typical in the functional setup and show up essentially due to the non-uniqueness of eigenfunctions but apparently do not affect the statistic \eqref{eq:ml_stat_F}. 
As already indicated, above assumptions are sufficient to obtain the limiting distribution of $\hat{T}_n$ which is stated in the next proposition and allows us to obtain critical values.
\begin{theorem}\label{thm:est_limit} Let Assumptions \eqref{eq:kernelL2}, (L), (P1), (A1), (A2) and $\lambda_d>0$ hold true. Then under $H_0$ it holds that, as $n\rightarrow \infty$,
\begin{equation}\label{eq:limit_tnhat}
	\lim_{n\rightarrow\infty} P\left(a(\log n)\hat{T}_n- b_d(\log n)\leq x\right)=\exp(-2\exp(-x))
\end{equation}
for all $x\in\mR$.
\end{theorem}
\begin{remark}
By considering the univariate analogue of \autoref{ass:assumptionA2} we note that $g(n)=(\log \log n)^{1/2}$ is the best possible rate (cf., e.g.\ \citet[Theorem A.4.1]{horvath1997limit} for the famous Darling-Erd{\H o}s asymptotics for i.i.d.\ random variables) and if such a rate holds true in \autoref{ass:assumptionP1}, then the rates in  \autoref{ass:assumptionA1} and \autoref{ass:assumptionA2} coincide and are both of order $o_P((\log\log n)^{-1})$.\footnote{Clearly, \autoref{ass:assumptionP1} with $g(n)=(\log \log n)^{1/2}$, could also be deduced from a law of the iterated logarithm for $\{Y_i\}_{i\in \mZ}$ and for the time-inversed counterpart $\{Y_{-i}\}_{i\in \mZ}$. However, to the best of our knowledge, the law of the iterated logarithm and results of Daling-Erd{\H o}s-type are not proven in the functional framework of \lpm-approximability so far.} 
In this article we will discuss results that allow us {\frqq a mathematically convenient derivation\flqq} of logarithmic rates $g(n)$ which are slightly weaker than $(\log \log n)^{1/2}$ but more than sufficient for our {\frqq practical\flqq} purposes (cf. \autoref{prop:moments2max}, \autoref{prop:twosided} and \autoref{lem:partialsum}). 
\end{remark}

Before proceeding with the verification of conditions (L) - (A2), we turn to the alternative and state the conditions which
ensure that the procedure detects changes  \begin{align} \Delta(t)=\mu_n(t)-\mu_1(t)\label{eq:delta}
\end{align} 
with $\Delta\neq0$ in $L^2[0,1]$ with probability tending to 1, as $n\rightarrow\infty$. Analogous to \autoref{ass:assumptionP1}, we need a bound on partial sums and conditions on estimates $\hat{\lambda}_j$ and $\hat{v}_j$. 
\begin{taggedassumption}{(P2)}\label{ass:assumptionP2} The weak law of large numbers holds true, i.e.\ it holds that $\|\sum_{i=1}^n Y_i\|=o_P(n)$ as $n\rightarrow \infty$. 
\end{taggedassumption}
Estimates $\hat{\lambda}_j$ appear in the denominator of \eqref{eq:stat_est_p2} and therefore need to be bounded. Recall, that we tacitly assume that all estimates are non-negative and in a decreasing order. 
\begin{taggedassumption}{(B1)}\label{ass:assumptionB1} Under $H_A$ it holds that $\hat{\lambda}_1=o_P(n/(\log \log n))$ as $n\rightarrow \infty$.
\end{taggedassumption}
\begin{taggedassumption}{(B2)}\label{ass:assumptionB2} Under $H_A$ it holds that, as $n\rightarrow \infty$, 
\begin{equation}\label{eq:cond_lim_space}
	|\int{\!\Delta(t)\hat{v}_r(t)}dt| \stackrel{P}{\longrightarrow} \xi,
\end{equation}
for some $1\leq r \leq d$ and some $\xi>0$. The change-function $\Delta(t)$ is defined in \eqref{eq:delta}.
\end{taggedassumption}

In order to obtain consistency for the change estimate $\hat{k}$ we state a condition that extends \citet[Theorem 2.8.1, disp. (2.8.7)]{horvath1997limit} to projection-based estimates in the functional framework.
\begin{taggedassumption}{(E1)}\label{ass:assumptionE1} Under $H_A$ it holds that, as $n\rightarrow \infty$, 
\begin{equation}\label{eq:condition_matrix}
	\frac{n}{g^2(n)}\left(\hat{\lambda}_d/\hat{\lambda}_1\right)\stackrel{P}{\longrightarrow} \infty,
\end{equation}
where the rate function $g(n)$ is the same as in \autoref{ass:assumptionP1}.
\end{taggedassumption} 
 Note that $|\hat{\lambda}_1/\hat{\lambda}_d|$ is the condition number of $\hat{\bSigma}$ with respect to the Euclidean norm. Hence, one possible interpretation is that \eqref{eq:condition_matrix} excludes ill-conditioned estimates.
\begin{remark}
Condition \eqref{eq:cond_lim_space} has an intuitive interpretation in terms of \eqref{eq:proj_innov}. Therefore, notice that the condition $\int\!\Delta(t) v_r(t)dt = 0$ for all $1\leq r\leq d$ is equivalent to $\bmu_1=\ldots =\bmu_n$, i.e.\ there would be {\frqq asymptotically\flqq} no change in the projected times series $\{\bX_i\}$. Hence, \eqref{eq:cond_lim_space} means that the change $\Delta$ has to be {\frqq asymptotically visible\flqq} in the projected time series $\{\hat{\bX}_i\}$ in \eqref{eq:proj_innov}.
\end{remark}

The above conditions are sufficient to state the following consistency results. 
\begin{theorem}\label{thm:functional_consistency_F} Let Assumptions \eqref{eq:kernelL2}, (L), (P2), (B1), (B2) and $\lambda_d>0$ hold true. Then under $H_A$ it holds that, as $n\rightarrow \infty$,
\begin{equation}\label{eq:consist_eq}
	(\log \log n)^{-1/2}\hat{T}_n\stackrel{P}{\longrightarrow} \infty.
\end{equation}

\end{theorem}

\begin{theorem}\label{thm:cp_est_consistency} Let Assumptions \eqref{eq:kernelL2}, (P1), (B2), (E1) and $\lambda_d>0$ hold true. Then under $H_A$ it holds that, as $n\rightarrow \infty$,
\begin{equation}\label{eq:consist_eq}
	\hat{k}(\hat{v},\hat{\lambda},d)/n\stackrel{P}{\longrightarrow} \theta.
\end{equation} 
\begin{remark}[Fully-functional estimates]\label{rem:ffunc}
An obvious drawback of projection-based approaches are the assumptions on the eigenstructure and on the visibility of projected changes. These assumptions can be {\frqq relaxed\flqq} in two steps:
\begin{enumerate}
\item First, note that \autoref{thm:cp_est_consistency} may be immediately restated for the estimate
\[
\hat{k}(\hat{v},1,d):=\operatorname*{arg\,max}_{1 \leq  k < n}w(k/n)\left(\hat{\eta}^2_{k,1}+\ldots +\hat{\eta}^2_{k,d}\right)^{1/2}
\]
without any assumptions on the eigenvalue estimates, without requiring $\lambda_d>0$ and such that \autoref{ass:assumptionE1} simplifies to $g^2(n)/n\rightarrow 0$.   However, we still need the visibility of the projected changes of \autoref{ass:assumptionB2}. 
\item  One possibility to avoid this issue is shown recently by \citet{aue2015} in a closely related context. They considered fully-functional estimates of the change point to overcome problems with \textit{\frqq high-frequency\flqq} changes. Those changes are  {\frqq difficult\flqq} to detect with the principal component approach since they require a large dimension $d$ to satisfy \autoref{ass:assumptionB2}. Working with large $d$'s in turn requires to estimate small eigenvalues $\lambda_d\approx 0$, i.e.\ the change point estimation becomes {\frqq unstable\flqq}.\footnote{In some sense this instability is reflected by \autoref{ass:assumptionE1}.}  Hence, it is worth mentioning that relying only on \autoref{ass:assumptionM} and on \autoref{ass:assumptionP1} with $g^2(n)/n\rightarrow 0$ it is straightforward to show \eqref{eq:consist_eq} for the fully-functional Darling-Erd\H{o}s type estimate 
\begin{align*}
	\hat{k}(\hat{v},1,\infty)&:=\operatorname*{arg\,max}_{1 \leq  k < n}w(k/n)\lim_{d\rightarrow\infty}\left(\hat{\eta}^2_{k,1}+\ldots +\hat{\eta}^2_{k,d}\right)^{1/2} \\
&= \operatorname*{arg\,max}_{1 \leq  k < n} w(k/n)\big\|n^{-1/2}\sum_{i=1}^{k} ( X_i - \bar{X}_{n})\big\|.
\end{align*}
This follows by going through the proof of \autoref{thm:cp_est_consistency} or of the underlying result in \citet[Theorem 2.8.1]{horvath1997limit}.
 \end{enumerate}
\end{remark}

\end{theorem} 
\autoref{thm:functional_consistency_F} and \autoref{thm:cp_est_consistency} rely on \autoref{ass:assumptionB2}. The latter is verified under $H_A$ typically via the relation
\begin{equation}\label{eq:consist_story}
	\Big||\int\! \Delta(t)\hat{v}_r(t)dt|-|\int\! \Delta(t)w_r(t)dt|\Big|=\cO_P\big(\|\hat{v}_r-c_rw_r\|\big)=o_P(1)
\end{equation}
on showing that $\|\hat{v}_r-c_rw_r\|=o_P(1)$ for appropriate $w_1,\ldots w_r\in L^2[0,1]$ (where $c_r\in\{0,1\}$ are random) together with $\int\! \Delta(t)w_r(t)dt \neq 0$ for some $1\leq  r\leq d$.
In the i.i.d.\ setting natural estimates $\hat{v}_r$ are given by the functional empirical principal components  (cf. \citet{berkes2009functional}) and it is shown that they converge (up to signs) to eigenfunctions $w_r$ of a contaminated operator (cf. also \citet{kirch2012functional} and \citet{torg2014}). Therefore, however, technical conditions on the eigenstructure of the contaminated operator together with the orthogonality condition $\int\! \Delta(t)w_r(t)dt \neq 0$ for some $1\leq  r\leq d$ have to be additionally imposed. In our setup, estimates $\hat{v}_r$ can always be chosen such that $\int \!\Delta(t)w_r(t)dt \neq 0$ (and therefore \autoref{ass:assumptionB2}) is fulfilled even with $r=1$ where $w_1=\Delta/\|\Delta\|$. This is stated in \autoref{prop:consistency_n} and has been observed by \citet{horvath2013testing} in a related context.
\\
\\
We turn to the verification of the conditions stated in Assumptions (L), (P1) and (P2) in case of \lpm-approximability as described in \autoref{sect:weakdep}. Conditions of \autoref{ass:assumptionA1} and \autoref{ass:assumptionA2} as well as of \autoref{ass:assumptionB1} and \autoref{ass:assumptionB2} concerning estimation will be verified separately in the next section.

\begin{theorem}\label{thm:null_lpm} Let $\{Y_i\}_{i\in \mZ}$ be causal and \lpm-approximable with $\kappa>2$ and a rate $\delta(m)=m^{-\gamma}$ for some $\gamma>2$. Then Assumptions \eqref{eq:kernelL2}, (L) and (P2) are fulfilled.
\end{theorem}

The next proposition applies the famous work of \citet{moricz1976moment} which, in combination with a result of \cite {tomacsa2006hajek}, will allow us to establish \autoref{ass:assumptionP1}. 
\begin{proposition}\label{prop:moments2max}Assume that for some constant $\kappa>2$, it holds that, as $n\rightarrow \infty$,
\begin{equation}\label{eq:econd}
	\left[E \big\|\sum_{i=1}^n Y_i\big\|^{\kappa}\right]^{1/\kappa} = \cO(n^{1/2}).
\end{equation}
Then, it holds that, as $n\rightarrow\infty$,
\[
	\max_{1\leq  k \leq n}k^{-1/2} \big\|\sum_{i=1}^k Y_i\big\|= \cO_P((\log n)^{1/\kappa}).
\]
\end{proposition}
Now, in order to verify \autoref{ass:assumptionP1}, it is sufficient to show \eqref{eq:econd} for $\{Y_i\}_{i\in \mZ}$ and $\{Y_{-i}\}_{i\in \mZ}$.  For \lpm-approximable and {\it{causal}} time series, \eqref{eq:econd} follows from \citet[Theorem 3.3]{berkes2013rice} given that $\kappa\in (2,3)$. In this article we verify \eqref{eq:econd} for \lpm-approximable {\it{noncausal}} time series based directly on \citet[Proposition 4]{berkes2011split}.
\begin{proposition}\label{prop:twosided} Let $\{Y_i\}_{i\in \mZ}$ be \lpm-approximable (not necessarily causal) with $\kappa\in [2,3)$ and $\sum_{m=1}^{\infty}\delta(m)<\infty$. Then, as $n\rightarrow\infty$, it holds that
\[
	\left[E \big\|\sum_{i=1}^n Y_{i}\big\|^{\kappa}\right]^{1/\kappa} = \cO(n^{1/2}).
\]
\end{proposition}
A combination of \autoref{prop:moments2max} and \autoref{prop:twosided} immediately yields the following result.
\begin{corollary}\label{lem:partialsum} Let $\{Y_i\}_{i\in \mZ}$ be \lpm-approximable with $\kappa>2$ and causal with $\sum_{m=1}^{\infty}\delta(m)<\infty$, then \autoref{ass:assumptionP1} is fulfilled with rate $g(n)=(\log n)^{1/2}$.
\end{corollary}

Altogether, \autoref{thm:null_lpm} and  \autoref{lem:partialsum} verify conditions \eqref{eq:kernelL2}, (L), (P1) and (P2) under \lpm-approximability with $\delta(m)=\cO(m^{-\gamma})$ for $\gamma>2$. The remaining Assumptions (A1), (A2), (B1) and (B2) ensure the validity of \eqref{eq:limit_tnhat} and of \eqref{eq:consist_eq}, in view of \autoref{thm:est_limit} and \autoref{thm:functional_consistency_F}. All these assumptions can be verified using suitable estimates which will be shown in the next section.

\section{Estimation of the eigenstructure}\label{sect:estimation}

In this section we discuss suitable estimates $\{\hat{v}_i\}$ and $\{\hat{\lambda}_i\}$ for $\{v_i\}$ and $\{\lambda_i\}$ which, as pointed out by \citet{reeder2012twoB}, is an intricate problem. One possibility to obtain such estimates, suggested by the latter, is to consider the eigenstructure of Bartlett-type estimators $\hat{\cC}$ for $\cC$ of the following general form
\[
 \big(\hat{\cC}x\big)(t)= \int{\!\hat{\zeta}(t,s) x(s)} ds,
 \]
where $x \in L^2[0,1]$, $t\in[0,1]$. The kernel $\hat{\zeta}$ is given by
\begin{equation}\label{eq:kernel_upn}
\hat{\zeta}(t,s) = \sum_{r=-n}^{n}K(r/h_n)\hat{\zeta}_r(t,s),\qquad t,s \in[0,1],
\end{equation}
with covariance estimators  
\[
 \hat{\zeta}_r(t,s)=
\begin{cases}\frac{1}{n}\sum_{i=1}^{n-r}  \left( X_{i}(t)-\bar{X}_n(t)\right) \left(X_{i+r}(s)-\bar{X}_n(s)\right),& r\geq 0,\\
\hat{\zeta}_{-r}(s,t),& r<0,
\end{cases}
\]
a symmetric, bounded and compactly supported kernel function $K(x)$ with $K(0)=1$ and a bandwidth $h_n \rightarrow \infty$ fulfilling $h_n=o(n)$. (Notice that, due to the compact support of $K(x)$, the summation in \eqref{eq:kernel_upn} is only up to $\lfloor ch_n\rfloor$ for some $c>0$.) These estimates were explored by \citet{reeder2012twoB} under \lpm-approximability in the context of a related two-sample problem.

We restrict ourselves now  to the framework of \lpm-approximable time series. For the sake of simplicity we consider exponential decay rates $\delta(m)=\exp(-cm)$, $c>0$. This is not very restrictive and already covers important time series models, in particular the functional AR($p$) time series, and will be sufficient for our purposes.   
Notice that $\hat{\cC}$ is symmetric Hilbert-Schmidt, hence has a spectral decomposition 
\[
	(\hat{\cC} x)(t) = \sum_{j=1}^{\infty}\hat{\lambda}_{j}\Big[\int\!x(s)\hat{v}_{j}(s)ds\Big] \hat{v}_{j}(t) 
\]
with real eigenvalues  $\hat{\lambda}_1\geq\hat{\lambda}_2\geq \ldots$ and corresponding orthonormal eigenfunctions $\hat{v}_1(t),\hat{v}_2(t),\ldots$.  Generally, the estimates of the eigenvalues may be negative but (at least under $H_0$) become eventually positive as $n\rightarrow\infty$.\footnote{A simple {\frqq ad-hoc\flqq} solution to ensure positiveness in finite samples (under $H_0$ and $H_A$) is to use the absolute values of these estimates  - this will be tacitly assumed in our subsequent theoretical considerations.} Subsequently, we use the same notation for the eigenstructure as used for generic estimates before. This should not lead to any confusion.

The following \autoref{thmp2:estimation_longrun} is an extension of Theorem 2 of \citet{reeder2012twoB}, where consistency has been shown (under weaker assumptions). We consider the case of \lfourm-approximable time series which allows us to work with the variances of the estimates.
\begin{theorem}\label{thmp2:estimation_longrun}Let $\{Y_i\}_{i\in \mZ}$ be \lfourm-approximable with $\delta(m)=\exp(-cm)$, $c>0$, and not necessarily causal. Assume that, as $x\rightarrow 0$,
\[
	|K(x)-1|=\cO(x^\rho)
\] 
for some $\rho\geq 1$. Then under $H_0$, it holds that, $\|\hat{\zeta}-\zeta\|=\cO_P((h_n/n)^{1/2}h_n  + 1/h_n^\rho)$, as $n\rightarrow\infty$, using a bandwidth $h_n =\lfloor c'n^{1/\gamma}\rfloor$ with some $\gamma>3$ and $c'>0$.
\end{theorem}
As one would expect, the rate of convergence in \autoref{thmp2:estimation_longrun} reflects that a higher smoothness of the kernel $K(x)$ at $x=0$ reduces the contribution of the {\frqq bias\flqq}. Admissible values for the \textit{\frqq smoothness parameter\flqq} $\rho$ are, e.g.,
\begin{itemize} 
\item $\rho=1$ for the Bartlett kernel $K(x)=(1-|x|)1_{[-1,1]}(x)$, 
\item $\rho=2$ for the Parzen kernel
\item and arbitrary large $\rho$ for the flat-top kernels. 
\end{itemize}
The main implication of this \autoref{thmp2:estimation_longrun} for us is that a polynomial rate $\|\hat{\zeta}-\zeta\|$$=\smash{O_P(n^{-\varepsilon})}$ holds true for some $\varepsilon>0$. This allows a convenient verification of \autoref{ass:assumptionA1} and \autoref{ass:assumptionA2} via \autoref{lem:partialsum} together with Lemmas 2.2 and 2.3 of \citet{horvath2012fda}. For \autoref{ass:assumptionA2} we have to assume additionally, as common in the functional setup, that the first $d$ eigenvalues of $\cC$ are simple, i.e.\ that
\begin{equation}\label{eq:eigenvalues_distinct}
	\lambda_1 > \lambda_2 > \ldots > \lambda_d>\lambda_{d+1}
\end{equation}
holds true. Hence, decomposition \eqref{eq:cov_series} is unique up to signs. 

\begin{proposition} Under $H_0$ and the assumptions of \autoref{thmp2:estimation_longrun} the \autoref{ass:assumptionA1} holds true. If additionally \eqref{eq:eigenvalues_distinct} is assumed, then \autoref{ass:assumptionA2} also holds true.
\end{proposition}
We conclude this section by an observation, which follows in view of \eqref{eq:consist_story} and due to the Lemma B.2 of \citet{horvath2013testing} (cf. also (3.5) and (3.6) therein). 

\begin{proposition}\label{prop:consistency_n} Let $\{Y_i\}_{i\in \mZ}$ be \lfourm-approximable and causal. Under $H_A$ and the assumptions of \autoref{thmp2:estimation_longrun} the \autoref{ass:assumptionB1} and \autoref{ass:assumptionB2} hold true.
\end{proposition}

As already mentioned in \autoref{rem:different_cond} the dependence condition in \citet[disp. (2.4)]{horvath2013testing} is slightly different. However, Lemma B.2 can be restated under our conditions due to \citet[Theorem 1.2]{jirak2013}.

\section{Simulations}\label{sect:simulAX}
\subsection{Monte Carlo simulation}

We proceed with a Monte Carlo simulation of the finite sample behavior. In order to describe our setting and implementation details we recall that $X_i(t)=\mu_i(t)+Y_i(t)$. 

\paragraph{Simulation setup:}The signal $\mu_i$ is set to $\mu_i(t)\equiv 0$ for $i=1,\ldots,n$ under the null hypothesis and 
\begin{equation}\label{eq:iteration_decomp}
\mu_i(t)= 
\begin{cases} 0,  &i=1,\ldots,\lfloor n/2\rfloor, \\ 
		\sin(t), &i=\lfloor n/2\rfloor+1,\ldots, n
\end{cases} 
\end{equation}
under the alternative. The innovations follow the formal functional AR($1$) model
\begin{equation}\label{eq:iteration}
	Y_{i}(t) = \int\! \Psi(t,s) Y_{i-1}(s)ds + \varepsilon_i(t),
\end{equation}
for $t\in [0,1]$ and $i\in\mZ$, where the shocks $\{\varepsilon_i\}$ are assumed to be Gaussian. Under the assumption of $\|\Psi\|<1$ this equation is known to have a unique \lpm-approximable solution where $\kappa\geq 2$ is arbitrary (due to Gaussianity of $\varepsilon_i$'s) and where the decay rate $\delta(m)$ is exponential according to \citet{hoermann2010weakly}.
 We will analyze the performance using two different kernels
\[
	\Psi_G(t,s)=C_G\exp((t^2+s^2)/2), \qquad \Psi_W(t,s)=C_W \min(t,s)
\]
normalized by constants $C_G,C_W\geq 0$, such that $\|\Psi\|=\psi$ for a prescribed value $\psi\in[0,1)$. These kernels are common benchmarks in the functional data change-point literature where $\Psi_G$ and $\Psi_W$ are usually referred to as \textit{\frqq Gaussian\flqq} or \textit{\frqq Wiener\flqq} kernels, respectively (cf. \citet{horvath2012fda}). 

\paragraph{Implementation details:} We have implemented the procedure in R using the {\frqq fda-package\flqq}. The shocks $\varepsilon_i(\cdot)$ are generated as paths of Brownian bridges on $[0,1]$ and are represented as functional objects via the fda-function  \verb+Data2fd(...)+ by using a \textit{\frqq$B$-Spline basis\flqq} of 25 functions. The same basis is also used to represent the kernel $\Psi$ and the innovations $Y_i(\cdot)$. More precise, the bivariate function $\Psi(t,s)$ is discretized on an equidistant grid $0=t_1<t_2<\ldots<t_T=1$ and for each $k=1,\ldots,T$ the univariate function $\Psi(t_k,\cdot)$ is then represented as a functional object. Next, 
\[
	I_i(t_k)=\int\! \Psi(t_k,s) Y_{i-1}(s)ds
\]
are computed for all $k=1,\ldots,T$ and $I_i(\cdot)$ itself is represented as a functional object with domain $[0,1]$, again, using the same $B$-Spline basis as in the previous steps. Having computed $I_i$ we are in the position to add up $Y_i(\cdot)=I_i(\cdot)+\varepsilon_i(\cdot)$. The correct representation of the dependence structure is ensured by using a so-called \textit{\frqq burn-in period\flqq} of length $N_B=100$, i.e.\ $Y_{1},\ldots Y_{n}$ are generated iteratively according to  \eqref{eq:iteration} beginning with $Y_{-{N_B}+1}:=\varepsilon_{-{N_B}+1}$ where, finally, the first $N_B$ observations $Y_{-{N_B}+1},\ldots, Y_{0}$ are discarded. The desired observations $X_{1},\ldots X_{n}$ are then created according to the signal plus noise representation $X_i=\mu_i + Y_i$ via \eqref{eq:iteration_decomp}.

The computation of the long run covariance estimator is carried out following \citet{reeder2012twoA} by using 25 orthonormal \textit{\frqq Fourier basis\flqq} functions. For simplicity we make use of a plain kernel $K(x)=1_{[-1,1]}(x)$, which yields satisfactory results. The overall picture remains comparable if one chooses e.g.\ a flat-top kernel as in \citet[disp. (4.1)]{reeder2012twoB} or a Bartlett kernel, instead.

\paragraph{Critical values:} The convergence in \eqref{eq:limit_tnhat} is rather slow. For that reason, we follow the idea investigated by \citet{horvath1997limit}, also successfully applied in a functional setting by \citet{torg2014}, by using quantiles of
\begin{equation}\label{eq:qvostp2}
V_n=\sup_{t\in I_n}\frac{|\bB^d(t)|}{(t(1-t))^{1/2}},
\end{equation}
where $\{\bB^d(t), 0\leq t\leq 1\}$ is a $d$-dimensional Gaussian process with components given by independent standard Brownian bridges $\{B_i(t),0\leq t\leq 1\}$ and e.g.\  $I_n=[h_n,1-h_n]$ with $h_n=(\log n)^{3/2}/n$. Asymptotic correctness of this choice follows from \eqref{eq:limit_tnhat} (cf. \citet[Corollary 1.3.1]{horvath1997limit} and the proof of \citet[Corollary 4.3]{torg2014}).\footnote{Our simulations confirm that critical values based on the Brownian-bridge type approximation \eqref{eq:qvostp2} outperform those based on the basic Gumbel-type limit \eqref{eq:limit_tnhat}. A (heuristic) reason is provided via Theorem 4.2 in \citet{torg2014} in a closely related setting. The counterpart of the latter theorem is beyond the scope of this article.} An essential advantage of \eqref{eq:qvostp2} is that quantiles can be computed using the expansion
\begin{equation}\label{eq:approx_vostrikova}
P\left(V_n \geq x\right)= \frac{x^d\exp(-x^2/2)}{2^{d/2}\Gamma(d/2)}\left\{\left(1-\frac{d}{x^2}\right)\log\frac{(1-h_n)^2}{h_n^2}+\frac{4}{x^2}+\cO(x^{-4})\right\}.
\end{equation}
This representation is well known as \textit{\frqq Vostrikova's tail approximation\flqq} (see \citet[disp. (18)]{vostrikova1981disorder} and also \citet[disp. (1.3.26)]{horvath1997limit}).

\paragraph{Dimension and bandwidth selection:}
Parameters $d$ and $h=h_n$ remain to be specified where especially the selection of $h$ is known to be a complex problem in practice. For example $d$ can be chosen according to the generalized CPV-Criterion (cf. Section 4.1 of \citet{horvath2013testing}) and $h$ could be specified (in appropriate cases) guided by rules from scalar time series as demonstrated by \citet{hoermann2010weakly}. However, both issues are not the focus of our research and therefore an overview for a range of parameters is presented in the tables below. 

\paragraph{Brief summary of simulations:} The behavior under $H_0$ or under $H_A$, respectively, is demonstrated in \autoref{Q} - \autoref{K} (based on 1000 repetitions). For moderate dependence and moderate sample sizes the procedure performs rather well and comparable for both kernels $\Psi_W$, $\Psi_G$. Moreover, it shows overall robustness with respect to the selection of various (small) dimensions $d$ and bandwidths $h$. With increasing sample size the {\frqq bias\flqq} due to the dependencies fades out which is in accordance with the nonparametric nature of the procedure.  
\begin{table}
\vspace{-1.5cm}
\centering
\caption{{Empirical sizes for functional AR(1) with Gaussian kernel $\Psi_G$; nominal level of $10\%$.}}
\begin{scriptsize}
\begin{tabular}{@{}lcrrrrrrrrrrrrr@{}} $n$ &  $\|\Psi_G\|$ & $\bf{d}$& \bf{1} & \bf{2} & \bf{3} & \bf{4} & \bf{5} & \phantom{d}& \bf{1} & \bf{2} & \bf{3} & \bf{4} & \bf{5}\\ \toprule && & \multicolumn{5}{c}{$h=1$}&& \multicolumn{5}{c}{$h=2$}\\\midrule
50 & 0.1 && 7.2& 6.9& 5.1& 3.4& 3.9& & 9.9& 6.7& 4.9& 5.1& 5.2\\ \midrule
& 0.2 && 8.7& 6.2& 4.3& 3.6& 2.6& & 8.2& 6.5& 5.9& 6.2& 5.8\\ \midrule
& 0.4 && 8.7& 5.3& 3.4& 3.9& 2.3& & 7.0& 5.8& 5.3& 5.3& 5.7\\ \midrule
& 0.6 && 8.9& 5.2& 3.2& 2.8& 2.3& & 4.7& 2.3& 3.1& 3.2& 4.0\\ \midrule
& 0.8 && 11.8& 7.9& 4.7& 2.9& 2.9& & 6.4& 3.1& 3.0& 3.1& 3.8\\ \midrule
\midrule
100 & 0.1 && 11.5& 10.1& 9.7& 8.6& 6.5& & 9.8& 9.4& 7.9& 6.1& 5.9\\ \midrule
& 0.2 && 7.8& 7.0& 5.5& 5.4& 4.6& & 9.0& 7.3& 6.0& 5.1& 4.8\\ \midrule
& 0.4 && 9.3& 6.8& 5.6& 5.9& 4.9& & 6.7& 4.6& 5.1& 3.8& 3.3\\ \midrule
& 0.6 && 14.6& 9.3& 8.8& 6.1& 4.8& & 6.7& 4.9& 3.9& 4.5& 4.0\\ \midrule
& 0.8 && 16.1& 11.1& 8.2& 7.0& 5.2& & 7.1& 6.0& 4.6& 4.1& 3.5\\ \midrule
\midrule
300 & 0.1 && 10.8& 10.5& 12.1& 10.8& 9.6& & 11.1& 9.3& 7.6& 6.9& 7.8\\ \midrule
& 0.2 && 9.6& 10.1& 8.9& 8.5& 8.2& & 9.0& 9.8& 8.0& 7.7& 7.0\\ \midrule
& 0.4 && 10.9& 10.0& 9.5& 9.6& 9.3& & 9.7& 10.9& 8.5& 7.8& 8.5\\ \midrule
& 0.6 && 14.4& 12.0& 9.1& 8.4& 8.2& & 7.4& 6.9& 6.4& 6.5& 6.7\\ \midrule
& 0.8 && 20.2& 15.2& 12.9& 12.8& 9.8& & 10.1& 8.8& 7.9& 8.1& 7.4\\ \midrule
\midrule
500 & 0.1 && 10.8& 9.3& 9.4& 9.0& 8.9& & 10.1& 8.8& 8.3& 7.9& 9.1\\ \midrule
& 0.2 && 9.7& 9.1& 9.0& 9.7& 9.6& & 8.1& 8.3& 6.8& 6.4& 7.1\\ \midrule
& 0.4 && 11.7& 11.2& 10.8& 10.2& 10.1& & 9.5& 9.5& 9.9& 9.0& 8.4\\ \midrule
& 0.6 && 15.5& 12.7& 11.1& 11.3& 10.2& & 10.3& 9.2& 8.3& 7.8& 7.8\\ \midrule
& 0.8 && 20.6& 18.1& 15.6& 14.1& 13.0& & 12.4& 9.1& 10.3& 9.1& 8.3\\
\toprule 
\\
 && & \multicolumn{5}{c}{$h=3$}&& \multicolumn{5}{c}{$h=4$}\\\midrule
50 & 0.1 && 10.0& 8.3& 7.2& 8.1& 8.5& & 9.7& 8.0& 9.7& 11.3& 11.9\\ \midrule
& 0.2 && 9.1& 6.3& 7.1& 7.2& 8.3& & 9.3& 7.3& 8.7& 9.1& 10.2\\ \midrule
& 0.4 && 6.5& 4.0& 5.1& 6.6& 8.4& & 7.2& 6.1& 6.6& 7.4& 10.1\\ \midrule
& 0.6 && 4.8& 3.4& 3.1& 3.4& 3.7& & 5.9& 5.3& 5.8& 6.8& 9.5\\ \midrule
& 0.8 && 2.7& 2.3& 2.3& 2.5& 4.9& & 4.3& 4.1& 5.0& 5.8& 8.3\\ \midrule
\midrule
100 & 0.1 && 9.8& 7.1& 6.3& 6.2& 6.1& & 9.2& 8.1& 6.6& 6.1& 7.0\\ \midrule
& 0.2 && 8.7& 7.9& 7.4& 7.6& 7.4& & 9.3& 8.4& 7.3& 7.5& 7.6\\ \midrule
& 0.4 && 8.2& 5.5& 5.0& 5.7& 5.9& & 6.7& 5.3& 6.7& 6.1& 7.7\\ \midrule
& 0.6 && 5.5& 3.7& 4.1& 3.9& 4.3& & 7.1& 5.7& 5.0& 4.5& 5.0\\ \midrule
& 0.8 && 5.9& 3.8& 3.6& 3.5& 5.6& & 4.1& 3.0& 3.4& 4.1& 4.8\\ \midrule
\midrule
300 & 0.1 && 9.5& 10.6& 9.7& 8.3& 8.6& & 10.6& 10.0& 10.6& 9.2& 7.5\\ \midrule
& 0.2 && 7.5& 10.0& 9.7& 9.0& 8.5& & 8.5& 7.9& 8.3& 8.2& 7.5\\ \midrule
& 0.4 && 8.5& 7.3& 6.6& 5.3& 5.4& & 9.1& 10.0& 9.1& 7.6& 6.7\\ \midrule
& 0.6 && 6.9& 6.8& 7.3& 6.5& 6.8& & 8.5& 7.1& 7.6& 6.5& 5.3\\ \midrule
& 0.8 && 6.7& 7.8& 7.7& 6.9& 7.1& & 6.7& 4.7& 7.4& 6.5& 6.1\\ \midrule
\midrule
500 & 0.1 && 11.6& 9.8& 7.7& 8.5& 8.8& & 10.3& 10.7& 9.3& 9.2& 8.1\\ \midrule
& 0.2 && 9.5& 8.7& 7.5& 8.1& 6.7& & 9.9& 9.7& 8.1& 7.9& 8.5\\ \midrule
& 0.4 && 7.9& 7.9& 7.0& 6.3& 7.7& & 7.2& 8.6& 8.6& 8.2& 7.9\\ \midrule
& 0.6 && 7.4& 6.7& 7.7& 7.5& 7.4& & 8.6& 7.1& 6.7& 6.8& 6.5\\ \midrule
& 0.8 && 7.2& 6.4& 7.2& 6.0& 7.1& & 7.3& 7.9& 7.7& 6.9& 5.9\\\bottomrule  
\end{tabular}
\end{scriptsize} 
\label{Q}
\end{table}

\begin{table} 
\vspace{-1.5cm}
\centering
\caption{{Empirical sizes for functional AR(1) with Wiener kernel $\Psi_W$; nominal level of $10\%$.}}
\begin{scriptsize}
\begin{tabular}{@{}lcrrrrrrrrrrrrr@{}} $n$ &  $\|\Psi_W\|$ & $\bf{d}$& \bf{1} & \bf{2} & \bf{3} & \bf{4} & \bf{5} & \phantom{d}& \bf{1} & \bf{2} & \bf{3} & \bf{4} & \bf{5}\\ \toprule && & \multicolumn{5}{c}{$h=1$}&& \multicolumn{5}{c}{$h=2$}\\\midrule
50 & 0.1 && 9.6& 7.6& 4.6& 3.9& 2.9& & 10.1& 8.0& 6.8& 6.4& 5.9\\ \midrule
& 0.2 && 9.8& 5.9& 4.6& 2.8& 2.8& & 8.0& 4.4& 4.0& 4.8& 5.1\\ \midrule
& 0.4 && 9.9& 5.7& 3.8& 2.9& 2.3& & 5.2& 3.5& 2.8& 2.4& 3.1\\ \midrule
& 0.6 && 11.2& 5.2& 2.7& 2.1& 1.8& & 4.8& 3.3& 2.5& 2.7& 3.2\\ \midrule
& 0.8 && 20.6& 11.0& 4.8& 3.0& 2.5& & 3.8& 1.9& 1.4& 2.1& 2.1\\ \midrule
\midrule
100 & 0.1 && 8.6& 8.1& 6.5& 6.7& 6.5& & 9.0& 8.7& 6.6& 7.1& 6.0\\ \midrule
& 0.2 && 8.8& 8.0& 7.9& 6.0& 5.8& & 8.5& 6.7& 6.1& 5.1& 4.9\\ \midrule
& 0.4 && 10.6& 8.6& 6.1& 5.1& 5.0& & 5.7& 7.1& 4.5& 4.4& 2.9\\ \midrule
& 0.6 && 13.9& 10.3& 7.4& 4.8& 3.9& & 6.0& 4.6& 3.2& 2.8& 2.7\\ \midrule
& 0.8 && 28.6& 19.6& 13.6& 9.7& 7.5& & 11.5& 7.0& 6.1& 4.1& 3.3\\ \midrule
\midrule
300 & 0.1 && 10.5& 9.4& 9.1& 8.7& 7.3& & 10.5& 9.9& 9.0& 6.8& 7.0\\ \midrule
& 0.2 && 9.6& 10.3& 10.0& 8.5& 8.9& & 9.2& 8.5& 8.3& 7.6& 7.2\\ \midrule
& 0.4 && 10.3& 9.8& 10.0& 8.5& 6.9& & 7.0& 6.9& 7.3& 5.7& 5.7\\ \midrule
& 0.6 && 19.1& 15.5& 12.7& 9.9& 8.5& & 8.9& 8.2& 7.7& 6.4& 6.0\\ \midrule
& 0.8 && 33.4& 25.5& 20.7& 18.6& 15.4& & 15.7& 12.5& 10.5& 9.7& 7.6\\ \midrule
\midrule
500 & 0.1 && 12.5& 10.2& 9.0& 8.6& 8.6& & 10.8& 10.5& 10.3& 9.3& 9.4\\ \midrule
& 0.2 && 10.8& 10.0& 8.9& 10.1& 8.4& & 9.5& 8.8& 8.3& 7.4& 6.9\\ \midrule
& 0.4 && 13.7& 13.3& 11.3& 9.9& 8.9& & 9.3& 8.5& 9.3& 8.4& 7.9\\ \midrule
& 0.6 && 18.2& 14.3& 12.8& 11.7& 10.9& & 11.4& 9.5& 7.9& 7.5& 8.6\\ \midrule
& 0.8 && 37.9& 29.0& 22.4& 19.1& 16.7& & 17.9& 12.7& 9.9& 9.9& 8.7\\
\toprule 
\\
&& & \multicolumn{5}{c}{$h=3$}&& \multicolumn{5}{c}{$h=4$}\\\midrule
50 & 0.1 && 10.5& 7.6& 8.2& 7.9& 7.9& & 8.4& 8.3& 8.9& 10.3& 12.3\\ \midrule
& 0.2 && 9.2& 7.8& 6.4& 6.4& 7.5& & 10.2& 8.0& 9.0& 9.9& 10.2\\ \midrule
& 0.4 && 5.5& 4.4& 3.9& 6.3& 5.7& & 6.0& 6.2& 6.0& 8.1& 9.2\\ \midrule
& 0.6 && 4.0& 2.8& 3.5& 4.1& 3.9& & 3.4& 4.2& 5.2& 5.0& 6.1\\ \midrule
& 0.8 && 1.9& 3.3& 3.1& 3.9& 3.4& & 2.4& 2.8& 3.1& 3.5& 5.4\\ \midrule
\midrule
100 & 0.1 && 9.3& 8.3& 7.7& 6.7& 7.1& & 10.9& 6.9& 7.5& 7.8& 9.3\\ \midrule
& 0.2 && 9.0& 5.5& 5.4& 4.9& 5.5& & 8.5& 6.5& 6.6& 6.2& 7.0\\ \midrule
& 0.4 && 7.3& 4.1& 4.0& 3.9& 4.5& & 7.2& 4.9& 5.3& 5.8& 6.3\\ \midrule
& 0.6 && 5.5& 4.3& 4.1& 4.1& 4.0& & 5.1& 3.4& 3.3& 4.2& 4.8\\ \midrule
& 0.8 && 5.8& 3.9& 3.1& 3.3& 3.2& & 2.8& 3.3& 2.9& 4.3& 4.7\\ \midrule
\midrule
300 & 0.1 && 9.7& 8.4& 7.4& 6.5& 6.2& & 9.2& 9.0& 7.5& 6.8& 6.8\\ \midrule
& 0.2 && 9.3& 9.1& 7.9& 6.8& 6.9& & 10.1& 9.1& 8.2& 7.1& 6.4\\ \midrule
& 0.4 && 7.7& 8.1& 5.7& 5.8& 5.0& & 7.2& 8.3& 6.7& 6.7& 7.0\\ \midrule
& 0.6 && 8.9& 7.0& 6.4& 6.6& 5.7& & 6.2& 5.4& 5.7& 4.8& 5.8\\ \midrule
& 0.8 && 9.7& 7.7& 7.1& 7.1& 6.1& & 6.2& 7.3& 4.7& 6.0& 5.5\\ \midrule
\midrule
500 & 0.1 && 10.8& 9.8& 7.7& 8.4& 8.0& & 9.6& 8.6& 9.7& 8.5& 7.5\\ \midrule
& 0.2 && 7.4& 8.4& 7.0& 6.8& 5.9& & 7.8& 8.6& 7.0& 8.4& 8.8\\ \midrule
& 0.4 && 9.8& 8.9& 7.2& 7.1& 6.5& & 6.3& 7.9& 7.6& 8.7& 8.0\\ \midrule
& 0.6 && 9.1& 8.1& 6.3& 5.7& 5.7& & 6.7& 6.9& 8.8& 7.4& 6.4\\ \midrule
& 0.8 && 10.7& 8.8& 7.3& 7.9& 7.4& & 8.5& 7.8& 7.3& 7.0& 7.6\\\bottomrule
\end{tabular}
\end{scriptsize}
\label{W}
\end{table}

\begin{table}
\vspace{-1.5cm}
\centering
\caption{{Empirical power for functional AR(1) with Gaussian kernel $\Psi_G$; $\mu_1\equiv 0$ and $\mu_n(t)=\sin(t)$; nominal level of $10\%$.}}
\begin{scriptsize}
\begin{tabular}{@{}lcrrrrrrrrrrrrr@{}} $n$ &  $\|\Psi_G\|$ & $\bf{d}$& \bf{1} & \bf{2} & \bf{3} & \bf{4} & \bf{5} & \phantom{d}& \bf{1} & \bf{2} & \bf{3} & \bf{4} & \bf{5}\\ \toprule && & \multicolumn{5}{c}{$h=1$}&& \multicolumn{5}{c}{$h=2$}\\\midrule
50 & 0.1 && 99.4& 99.9& 99.9& 99.1& 52.4& & 95.2& 40.8& 3.3& 4.0& 4.0\\ \midrule
& 0.2 && 99.3& 99.8& 99.8& 98.6& 52.1& & 91.8& 37.5& 2.0& 2.2& 3.0\\ \midrule
& 0.4 && 94.8& 99.7& 99.1& 96.2& 45.0& & 83.0& 28.5& 1.5& 1.8& 2.8\\ \midrule
& 0.6 && 89.3& 98.7& 98.9& 94.7& 40.6& & 68.5& 23.0& 1.8& 1.7& 2.2\\ \midrule
& 0.8 && 81.1& 97.6& 96.1& 91.2& 38.7& & 54.0& 18.9& 0.6& 0.7& 1.7\\ \midrule
\midrule
100 & 0.1 && 100& 100& 100& 100& 100& & 100& 100& 100& 100& 100\\ \midrule
& 0.2 && 100& 100& 100& 100& 100& & 100& 100& 100& 100& 100\\ \midrule
& 0.4 && 100& 100& 100& 100& 100& & 100& 100& 100& 100& 100\\ \midrule
& 0.6 && 99.7& 100& 100& 100& 100& & 99.5& 100& 100& 100& 100\\ \midrule
& 0.8 && 98.7& 100& 100& 100& 100& & 95.4& 100& 100& 100& 100\\ \midrule
\midrule
300 & 0.1 && 100& 100& 100& 100& 100& & 100& 100& 100& 100& 100\\ \midrule
& 0.2 && 100& 100& 100& 100& 100& & 100& 100& 100& 100& 100\\ \midrule
& 0.4 && 100& 100& 100& 100& 100& & 100& 100& 100& 100& 100\\ \midrule
& 0.6 && 100& 100& 100& 100& 100& & 100& 100& 100& 100& 100\\ \midrule
& 0.8 && 100& 100& 100& 100& 100& & 100& 100& 100& 100& 100\\ \midrule
\midrule
500 & 0.1 && 100& 100& 100& 100& 100& & 100& 100& 100& 100& 100\\ \midrule
& 0.2 && 100& 100& 100& 100& 100& & 100& 100& 100& 100& 100\\ \midrule
& 0.4 && 100& 100& 100& 100& 100& & 100& 100& 100& 100& 100\\ \midrule
& 0.6 && 100& 100& 100& 100& 100& & 100& 100& 100& 100& 100\\ \midrule
& 0.8 && 100& 100& 100& 100& 100& & 100& 100& 100& 100& 100\\
\toprule 
\\  && & \multicolumn{5}{c}{$h=3$}&& \multicolumn{5}{c}{$h=4$}\\\midrule
50 & 0.1 && 45.3& 1.9& 2.8& 4.2& 5.3& & 1.1& 2.2& 2.7& 4.5& 6.3\\ \midrule
& 0.2 && 37.4& 2.1& 2.7& 3.7& 4.3& & 1.2& 2.3& 4.0& 4.4& 6.5\\ \midrule
& 0.4 && 28.2& 1.9& 2.8& 3.8& 5.3& & 1.1& 2.7& 3.9& 4.2& 7.2\\ \midrule
& 0.6 && 15.8& 1.0& 2.0& 3.0& 3.7& & 0.4& 1.6& 2.4& 4.3& 5.5\\ \midrule
& 0.8 && 9.7& 0.7& 1.0& 2.0& 3.4& & 1.0& 1.2& 2.9& 4.3& 5.8\\ \midrule
\midrule
100 & 0.1 && 100& 100& 100& 13.6& 3.8& & 99.9& 93.2& 2.2& 3.0& 2.8\\ \midrule
& 0.2 && 100& 100& 99.1& 12.7& 4.6& & 99.9& 91.6& 2.3& 3.0& 3.7\\ \midrule
& 0.4 && 99.9& 100& 98.9& 9.3& 2.9& & 99.3& 86.5& 2.5& 3.4& 3.9\\ \midrule
& 0.6 && 99.0& 100& 97.6& 8.9& 1.9& & 93.8& 81.0& 2.8& 3.9& 4.1\\ \midrule
& 0.8 && 91.5& 99.7& 95.8& 8.9& 2.6& & 82.1& 72.3& 0.6& 2.2& 2.2\\ \midrule
\midrule
300 & 0.1 && 100& 100& 100& 100& 100& & 100& 100& 100& 100& 100\\ \midrule
& 0.2 && 100& 100& 100& 100& 100& & 100& 100& 100& 100& 100\\ \midrule
& 0.4 && 100& 100& 100& 100& 100& & 100& 100& 100& 100& 100\\ \midrule
& 0.6 && 100& 100& 100& 100& 100& & 100& 100& 100& 100& 100\\ \midrule
& 0.8 && 100& 100& 100& 100& 100& & 100& 100& 100& 100& 100\\ \midrule
\midrule
500 & 0.1 && 100& 100& 100& 100& 100& & 100& 100& 100& 100& 100\\ \midrule
& 0.2 && 100& 100& 100& 100& 100& & 100& 100& 100& 100& 100\\ \midrule
& 0.4 && 100& 100& 100& 100& 100& & 100& 100& 100& 100& 100\\ \midrule
& 0.6 && 100& 100& 100& 100& 100& & 100& 100& 100& 100& 100\\ \midrule
& 0.8 && 100& 100& 100& 100& 100& & 100& 100& 100& 100& 100\\\bottomrule 
\end{tabular}
\end{scriptsize}
\label{Z}
\end{table}

\begin{table}
\vspace{-1.5cm}
\centering
\caption{{Empirical power for functional AR(1) with Wiener kernel $\Psi_W$; $\mu_1\equiv 0$ and $\mu_n(t)=\sin(t)$; nominal level of $10\%$.}}
\begin{scriptsize}
\begin{tabular}{@{}lcrrrrrrrrrrrrr@{}} $n$ &  $\|\Psi_W\|$ & $\bf{d}$& \bf{1} & \bf{2} & \bf{3} & \bf{4} & \bf{5} & \phantom{d}& \bf{1} & \bf{2} & \bf{3} & \bf{4} & \bf{5}\\ \toprule && & \multicolumn{5}{c}{$h=1$}&& \multicolumn{5}{c}{$h=2$}\\\midrule
50 & 0.1 && 99.0& 99.9& 99.7& 98.2& 54.3& & 93.3& 36.6& 2.2& 3.5& 4.1\\ \midrule
& 0.2 && 97.9& 99.7& 99.3& 97.2& 53.1& & 91.3& 31.4& 2.3& 3.1& 3.4\\ \midrule
& 0.4 && 93.0& 98.8& 99.0& 95.4& 43.4& & 73.7& 22.4& 2.0& 2.7& 2.7\\ \midrule
& 0.6 && 81.1& 93.4& 94.7& 90.9& 36.4& & 57.1& 12.9& 1.4& 1.2& 1.5\\ \midrule
& 0.8 && 73.2& 81.9& 89.1& 84.6& 40.3& & 37.3& 9.5& 1.5& 1.3& 1.6\\ \midrule
\midrule
100 & 0.1 && 100& 100& 100& 100& 100& & 100& 100& 100& 100& 100\\ \midrule
& 0.2 && 100& 100& 100& 100& 100& & 100& 100& 100& 100& 100\\ \midrule
& 0.4 && 99.9& 100& 100& 100& 100& & 99.6& 100& 100& 100& 100\\ \midrule
& 0.6 && 98.3& 100& 100& 100& 100& & 96.9& 100& 100& 100& 100\\ \midrule
& 0.8 && 93.4& 99.7& 100& 100& 100& & 85.2& 98.0& 99.9& 100& 100\\ \midrule
\midrule
300 & 0.1 && 100& 100& 100& 100& 100& & 100& 100& 100& 100& 100\\ \midrule
& 0.2 && 100& 100& 100& 100& 100& & 100& 100& 100& 100& 100\\ \midrule
& 0.4 && 100& 100& 100& 100& 100& & 100& 100& 100& 100& 100\\ \midrule
& 0.6 && 100& 100& 100& 100& 100& & 100& 100& 100& 100& 100\\ \midrule
& 0.8 && 100& 100& 100& 100& 100& & 100& 100& 100& 100& 100\\ \midrule
\midrule
500 & 0.1 && 100& 100& 100& 100& 100& & 100& 100& 100& 100& 100\\ \midrule
& 0.2 && 100& 100& 100& 100& 100& & 100& 100& 100& 100& 100\\ \midrule
& 0.4 && 100& 100& 100& 100& 100& & 100& 100& 100& 100& 100\\ \midrule
& 0.6 && 100& 100& 100& 100& 100& & 100& 100& 100& 100& 100\\ \midrule
& 0.8 && 100& 100& 100& 100& 100& & 100& 100& 100& 100& 100\\
\toprule 
\\
&& & \multicolumn{5}{c}{$h=3$}&& \multicolumn{5}{c}{$h=4$}\\\midrule
50 & 0.1 && 40.1& 1.7& 2.6& 3.1& 4.6& & 1.2& 2.1& 4.1& 4.6& 7.0\\ \midrule
& 0.2 && 31.6& 1.5& 2.2& 3.2& 4.4& & 0.8& 1.7& 2.1& 3.2& 6.0\\ \midrule
& 0.4 && 21.0& 2.1& 2.4& 2.8& 3.8& & 0.8& 2.1& 3.6& 5.0& 7.8\\ \midrule
& 0.6 && 10.5& 1.5& 2.2& 2.2& 3.3& & 0.5& 1.6& 1.6& 2.3& 5.0\\ \midrule
& 0.8 && 2.8& 0.5& 1.1& 1.6& 2.5& & 0.9& 1.0& 1.4& 1.6& 3.0\\ \midrule
\midrule
100 & 0.1 && 100& 100& 99.8& 10.5& 3.5& & 99.9& 92.5& 3.0& 3.6& 4.8\\ \midrule
& 0.2 && 99.9& 100& 99.4& 10.5& 3.5& & 99.7& 90.0& 2.5& 3.4& 3.3\\ \midrule
& 0.4 && 99.3& 99.7& 97.7& 10.9& 3.2& & 97.2& 76.7& 2.9& 2.6& 2.8\\ \midrule
& 0.6 && 93.0& 98.1& 95.4& 9.0& 2.9& & 81.4& 56.3& 1.7& 1.8& 2.9\\ \midrule
& 0.8 && 69.3& 92.4& 88.0& 9.2& 2.4& & 50.1& 46.3& 1.6& 1.7& 2.0\\ \midrule
\midrule
300 & 0.1 && 100& 100& 100& 100& 100& & 100& 100& 100& 100& 100\\ \midrule
& 0.2 && 100& 100& 100& 100& 100& & 100& 100& 100& 100& 100\\ \midrule
& 0.4 && 100& 100& 100& 100& 100& & 100& 100& 100& 100& 100\\ \midrule
& 0.6 && 100& 100& 100& 100& 100& & 100& 100& 100& 100& 100\\ \midrule
& 0.8 && 99.9& 100& 100& 100& 100& & 99.8& 100& 100& 100& 100\\ \midrule
\midrule
500 & 0.1 && 100& 100& 100& 100& 100& & 100& 100& 100& 100& 100\\ \midrule
& 0.2 && 100& 100& 100& 100& 100& & 100& 100& 100& 100& 100\\ \midrule
& 0.4 && 100& 100& 100& 100& 100& & 100& 100& 100& 100& 100\\ \midrule
& 0.6 && 100& 100& 100& 100& 100& & 100& 100& 100& 100& 100\\ \midrule
& 0.8 && 100& 100& 100& 100& 100& & 100& 100& 100& 100& 100\\\bottomrule 
\end{tabular}
\end{scriptsize}
\label{K}
\end{table}

\subsection{Application to Load Profiles}\label{sect:load_profile}
As a {\frqq real-life\flqq} example we take a closer look at electricity consumption data. This is inspired by the analysis of electricity data of \citet{rice2014}. We consider load profiles for the low voltage electricity network of the German electricity distribution company E.ON Mitte AG (now {\frqq EnergieNetz Mitte\flqq}) for 2012. 

\begin{figure}[H]
\centering
\includegraphics[width=\textwidth]{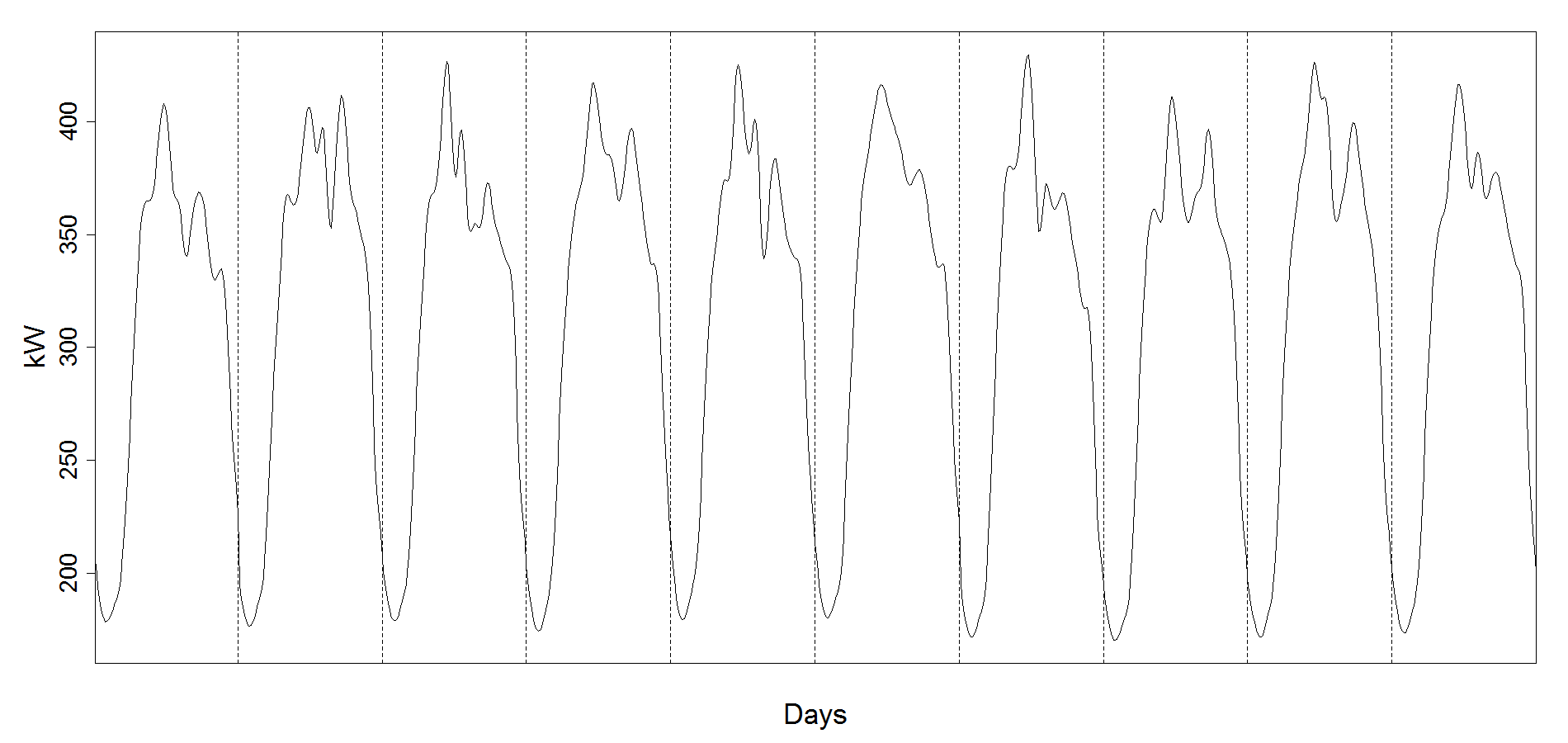}
\caption{10 consecutive daily load curves $X_{115}(t),\ldots,X_{124}(t)$ without weekends.}
\label{fig:daily_curve}
\end{figure}

\begin{figure}[H]
\centering
\includegraphics[width=\textwidth]{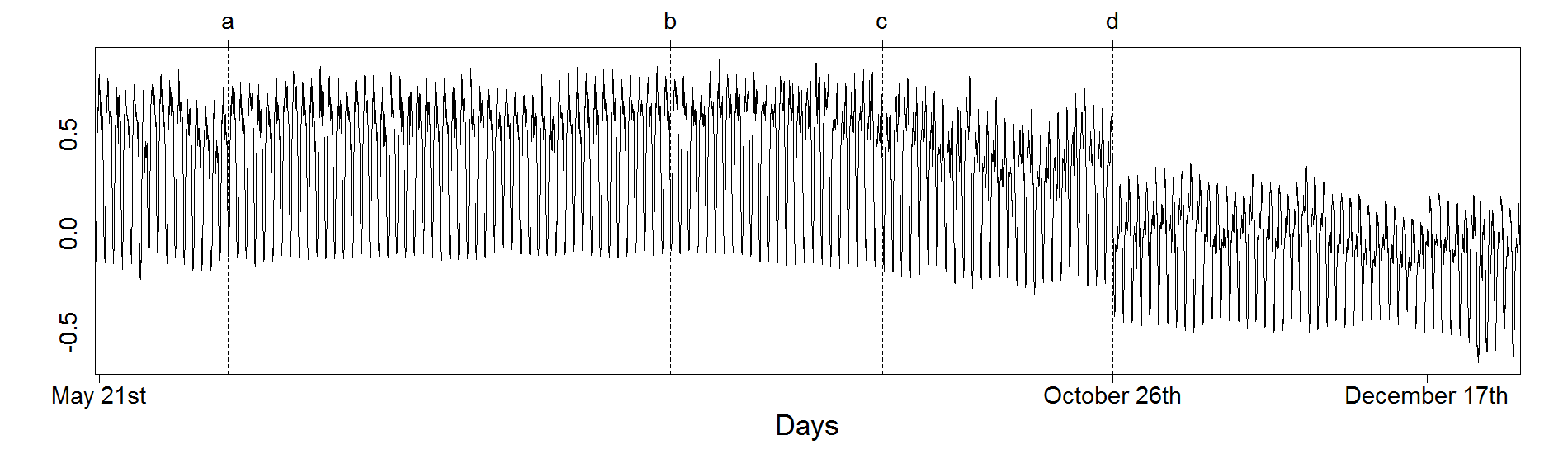}
\caption{Consecutive daily $\log$-transformed load curves $\tilde{X}_{101}(t),\ldots,\tilde{X}_{261}(t)$ corresponding to workdays in the time period May 21st - December 31st in 2012.}
\label{fig:profile_year_sequential}
\end{figure}

\paragraph{Data description:} The load profiles are based on quarter-hourly measurements (in kW), i.e.\ each of 366 days consists of 96 highly correlated observations. We split the original time series into segments corresponding to days and view each daily record as a curve, that is treat it as functional data (cf. \autoref{fig:daily_curve}). 
In order to apply the testing procedure we proceed as described in \autoref{sect:simulAX} and represent the discrete daily records as functional objects in rescaled time $t\in [0,1]$ using 25 $B$-Spline basis-functions, hereby smoothing the data. In the next step, in view of the obviously different stochastic pattern, we remove all curves which correspond to weekends. The remaining dataset consists of 261 curves $X_1(t),\ldots, X_{261}(t)$ corresponding to workdays (cf. \autoref{fig:profile_year_sequential}). Now, to gain stationarity, all curves are $\log$-transformed via 
\[
	\tilde{X}_i(t)=\log(X_i(t)/X_i(0)).
\] 
For a discussion on this transformation we refer to \citet{rice2014}. In a third step, we discard the observations $\tilde{X}_1(t),\ldots, \tilde{X}_{100}$ (i.e.\ the data before May 21st) which show a somewhat too erratic behavior to be reasonable for our analysis. 
\begin{figure}[H]
\centering
\includegraphics[width=\textwidth]{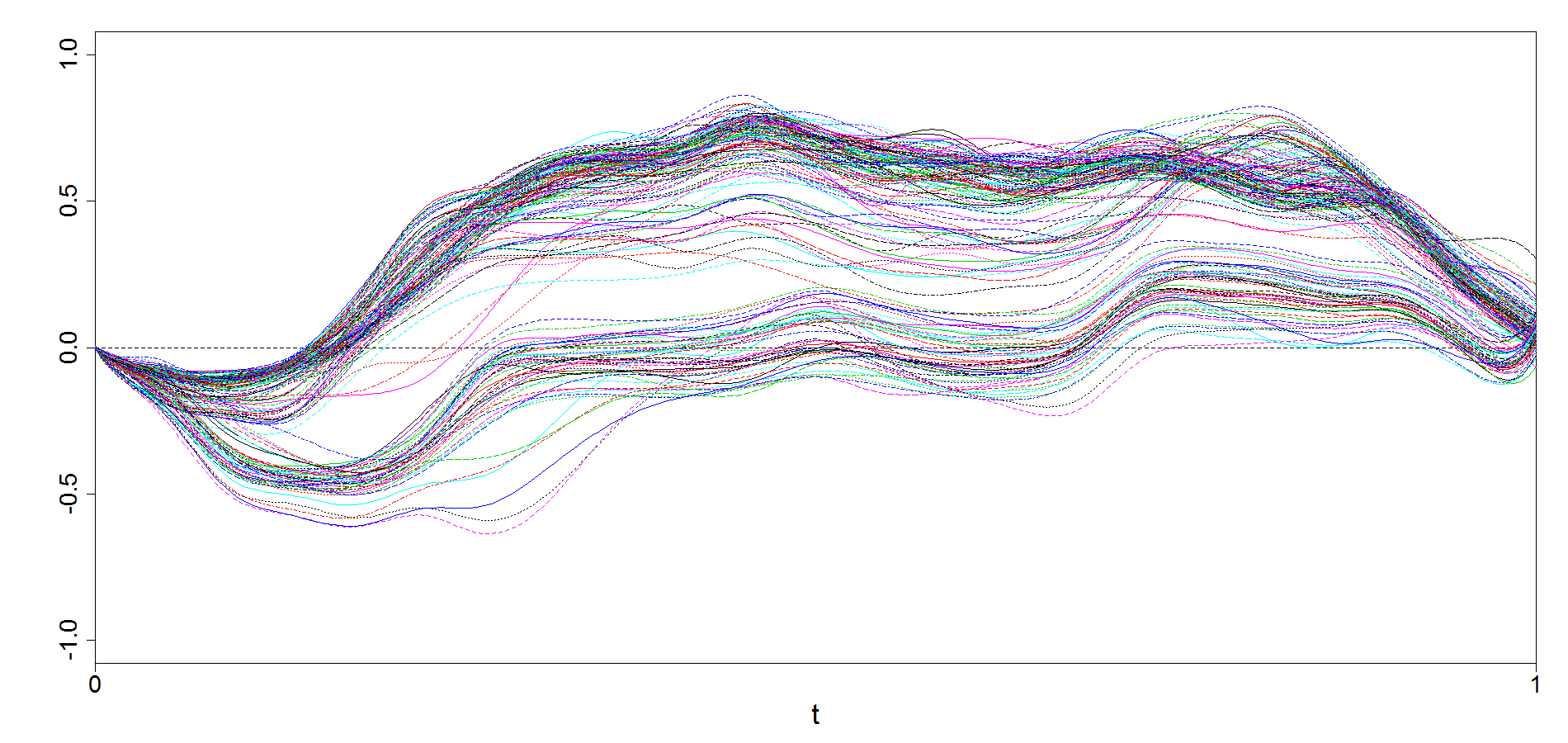}
\caption{Daily $\log$-transformed load curves $\tilde{X}_i(t)$ for the low voltage electricity grid corresponding to workdays in the time period May 21st - December 31st in 2012. Lower curves correspond to the winter months. This is the {\frqq functional representation\flqq} of the time series from \autoref{fig:profile_year_sequential}.}
\label{fig:profile_allinone}
\end{figure}
The remaining observations $\tilde{X}_{101}(t),\ldots,\tilde{X}_{261}(t)$ exhibit an obvious large abrupt change in the mean at $\tilde{X}_{216}$, i.e.\ at line {\frqq d\flqq} in \autoref{fig:profile_year_sequential}, and a quite stationary behavior before and after the jump.

 \paragraph{Data analysis:} Due to the large jump, we expect the procedure to reject the null hypothesis distinctly (the slight trend in segment {\frqq c - d\flqq} of \autoref{fig:profile_year_sequential} should not affect the performance), which is confirmed in \autoref{tab:sizeStatLoadProfileAX}: the procedure rejects the null hypothesis for a wide range of parameters (where the $p$-values are based on approximation \eqref{eq:approx_vostrikova}). The tests are carried out using a plain kernel $K(x)=1_{[-1,1]}(x)$, different bandwidths $h=0,\ldots,4$ and various subspace dimensions $d=1,\ldots,6$. (Recall that we divide by $h$ in \eqref{eq:kernel_upn} and $h=0$ is thus formally prohibited. Here, we set $\hat{\zeta}(t,s)=\hat{\zeta}_0(t,s)$ for $h=0$.) Note that $\hat{T}_n$ is largest (or vice versa the $p$-values are smallest) for $h=0$, i.e.\ when the dependence structure is not taken into account. However, the results for $h>0$ should be more reliable since there seem to be indicators of dependencies in the data described in the following. We performed a basic analysis and checked for independence using the functional \textit{\frqq Portmanteau Test of Independence\flqq} of \citet{gabrys2007} which is also based on a dimension-reduction approach using (static) functional principal components. First, note that this procedure already requires mean zero data. Therefore, to minimize the influence of obvious large changes and of less obvious smaller trends we restrict our considerations to the rather homogeneous segment $\tilde{X}_{115}(t),\ldots,\tilde{X}_{190}(t)$ (i.e.\ segment {\frqq a - c\flqq} in \autoref{fig:profile_year_sequential}) and center this subsample by its sample mean. The test for this sample yields small $p$-values $<10^{-10}$ for a range of parameters $d$ (number of principal components) and $\tilde{H}$ (maximum lag). We obtain somewhat larger, but still small, $p$-values (cf. \autoref{tab:TestIndependence}) if we restrict ourselves further to the segment $\tilde{X}_{115}(t),\ldots,\tilde{X}_{165}(t)$, (i.e.\ segment {\frqq a - b\flqq} in \autoref{fig:profile_year_sequential}.)

\begin{table}[H]
\caption{{Values of $\hat{T}_n$ for the load profile dataset. $p$-values are given in brackets.}}
\centering	 
\label{tab:sizeStatLoadProfileAX}
\begin{footnotesize}
\begin{tabular}{@{}lllllll@{}} 
\toprule 
$d$ 	&  	& 	$h=0$		&  $h=1$ &  $h=2$ &  $h=3$ &  $h=4$\\ 
\toprule
1 		& 	& 	 $12.09\;(<0.0001)$	&	$7.08\; (<0.0001)$		&	$5.53\; (<0.0001)$ &	$4.70\; (0.0001)$		&	$4.17\; (0.0014)$ \\
2 		& 	& 	 $12.17\;	 (<0.0001)$	&	$7.14\; (<0.0001)$	&	$5.58\; (<0.0001)$		&	$4.75\; (0.0007)$		&	$4.23\; (0.0057)$	\\
3		& 	& 	 $12.18\;	 (<0.0001)$	&	$7.14\; (<0.0001)$	&	$5.80\; (<0.0001)$		&	$4.75\; (0.0025)$		&	$4.23\; (0.0179)$	\\
4		& 	& 	 $12.34\;	 (<0.0001)$	&	$7.22\; (<0.0001)$		&	$5.91\; (<0.0001)$		&	$5.37\; (0.0005)$		&	$5.09\; (0.0018)$	\\
5 		& 	& 	$12.38	\; (<0.0001)$	&	$7.43\; (<0.0001)$		&	$5.92\; (<0.0001)$		&	$5.43\; (0.0011)$		&	$5.20\; (0.0029)$	 \\
6 		& 	& 	 $12.62\;	 (<0.0001)$	&	$7.51\; (<0.0001)$		&	$5.96\; (0.0002)$		&	$5.48\; (0.0021)$		&	$5.29\; (0.0048)$	\\
\bottomrule
\end{tabular}
\end{footnotesize}
\end{table}

\begin{table}[H] 
	\centering
\caption{{$p$-values from the Portmanteau Test of Independence of \citet{gabrys2007} applied to the segment $\tilde{X}_{115}(t),\ldots,\tilde{X}_{165}(t)$ of the load profile dataset. $d$ represents the number of principal components and $\tilde{H}$ denotes the maximum lag used for the test.}}
\label{tab:TestIndependence}
\begin{footnotesize}
\begin{tabular}{@{}llrrrr@{}} 
\toprule 
$\tilde{H}$ 	&  	& 	$d=1$		&  $d=2$ &  $d=3$ &  $d=4$\\
\toprule
1 		& 	& 	 $0.0011$		&	 $0.0228$		&	 	$0.0002$ 		&	 $0.0023$	\\
2 		& 	& 	 $<0.0001$ 	&	 $0.0025$		&	  	$<0.0001$		&	 $<0.0001$	\\
3		& 	& 	 $<0.0001$		&	 $0.0021$		&	 	$<0.0001$		&	 $<0.0001$\\
\bottomrule
\end{tabular}
\end{footnotesize}
\end{table}

\begin{remark}
Note that the curves in \autoref{fig:profile_year_sequential} and \autoref{fig:profile_allinone} which correspond to the winter months are below those corresponding to summer moths, due to the {\frqq functional rescaling\flqq} $X_i(t)/X_i(0)$ for $t\in[0,1]$: The electricity demand during morning and daytime in winter and summer months is comparable. However, in the winter the demand in the evening and especially at midnight, i.e.\ at $X_i(0)$, is much higher. Hence, the observed change is in accordance with the fact that electricity consumption in the winter is higher than in the summer and is most likely due to a switch in the {\frqq supply regime\flqq}. 
\end{remark}

\section{Proofs}\label{sect:proofII}

\begin{proof}[Proof of \autoref{thm:est_limit}]
We outline the important steps, thereby following the proof of \citet[Theorem 4.1]{torg2014}, which in turn is largely based on considerations of \citet{berkes2009functional}.

First, we replace the eigenfunctions with their estimates. Going through the proofs of Lemmas 6.4 and 6.5 of \citet{torg2014} we see that \autoref{ass:assumptionP1} and \autoref{ass:assumptionA2} ensure that, as $n\rightarrow \infty$,
\begin{equation}\label{eq:firstlim}
	|T_n(X;\hat{v},\lambda) - T_n(Y;v,\lambda)|=o_P((\log \log n)^{-1/2})
\end{equation}
and therefore, taking \autoref{ass:assumptionL} into account
\begin{equation}\label{eq:seclim}
	\lim_{n\rightarrow\infty} P\left(a(\log n)T_n(X;\hat{v},\lambda)- b_d(\log n)\leq x\right)=\exp(-2\exp(-x))
\end{equation}
holds true. Next, we replace the population eigenvalues with their empirical versions. \autoref{ass:assumptionA1} implies that $\lim_{n\rightarrow\infty}P(\hat{\lambda}_d>c)=1$, for some $c>0$, and that $\max_{i=1,\ldots,d}|\hat{\lambda}^{-1/2}_i-\lambda^{-1/2}_i|=o_P((\log\log n)^{-1})$ hold true. Following the arguments of \citet[Lemma 6.7]{torg2014} we see that \eqref{eq:seclim} and \autoref{ass:assumptionA1} imply that 
\[
	|T_n(X;\hat{v},\hat{\lambda}) - T_n(X;\hat{v},\lambda)|=o_P((\log \log n)^{-1/2}).
\]
The assertion follows immediately by using \eqref{eq:firstlim} and \autoref{ass:assumptionL}.
\end{proof}

\begin{proof}[Proof of \autoref{thm:functional_consistency_F}]
Let $m=\ntheta$. It is clear that
\begin{equation}\label{eq:decomp_ch}
	\sum_{i=1}^{k}\Big(X_i(t)-\bar{X}_n(t)\Big) = \sum_{i=1}^{k}\Big(Y_i(t)-\bar{Y}_n(t)\Big) - \frac{k}{n} (n-m)\Delta(t)
\end{equation}
holds true for $1\leq k \leq m$. Hence,  using standard arguments we obtain that
\begin{align*}
\hat{T}_n &\geq  c_1  w(m/n)\bigg\{\Big|\int \!n^{-1/2}\sum_{i=1}^{m} ( X_i(t) - \bar{X}_{n}(t))\hat{v}_r(t)dt \Big|\bigg\}/\hat{\lambda}^{1/2}_r\\
&\geq  c_2 w(\theta)\bigg\{\bigg| \big|\int \!n^{-1}\sum_{i=1}^{m} ( Y_i(t) - \bar{Y}_{n}(t))\hat{v}_r(t)dt \big|\\
&- \bigg[|\int\!\Delta(t)\hat{v}_r(t)dt|-\xi\bigg] m(n-m)/\big(n^{3/2}n^{1/2}\big)\\
&\quad +\xi m(n-m)/n^2\bigg|\bigg\}\times \bigg[(\log \log n)^{1/2} n^{1/2}/\big((\log \log n)^{1/2}\hat{\lambda}_r^{1/2}\big)\bigg] \\
&=: c_2 w(\theta) \big|A_1 + A_2 + A_3\big|(\log \log n)^{1/2} \big((\hat{\lambda}_r  \log \log n)/n\big)^{-1/2}
\end{align*}
for some $c_1,c_2>0$. We get $A_1=o_P(1)$ in view of \autoref{ass:assumptionP2}, $\|\hat{v}_r\|=1$ and due to the Cauchy-Schwarz inequality. Further, $A_2=o_P(1)$ holds true on account of \autoref{ass:assumptionB2}. The third term $A_3$ converges towards a nonzero (positive) constant, again due to $\xi>0$ in \autoref{ass:assumptionB2}. The assertion follows now in view of \autoref{ass:assumptionB1}.
\end{proof}

\begin{proof}[Proof of \autoref{thm:cp_est_consistency}]
We carry over and adapt the arguments of \citet[Theorem 2.8.1]{horvath1997limit} to our functional setting. \autoref{ass:assumptionE1} particularly implies that $P(\hat{\lambda}_d>0)$ is tending to 1 under $H_A$. Thus, for convenience we tacitly restrict the consideration to the set where $\hat{\lambda}_d>0$ holds true.

Following the notation in \eqref{eq:proj_innov} we write $\hat{\bX}_i$, $\hat{\bY}_i$ and $\hat{\bmu}_i$ for vectors where the $r$-th components are the scores $\hat{\bX}_{i,r}=\int\!X_i(t) \hat{v}_r(t)dt$,  $\hat{\bY}_{i,r}=\int\!Y_i(t) \hat{v}_r(t)dt$ and $\hat{\bmu}_{i,r}=\int\!\mu_i(t) \hat{v}_r(t)dt$. Furthermore, we set $\hat{\bDelta}=\hat{\bmu}_n-\hat{\bmu}_1$.  
Using the Cauchy-Schwarz inequality and \autoref{ass:assumptionP1} we have
\begin{align*}
	&\max_{1 \leq  k < n} w(k/n) \bigg|\int\! n^{-1/2}\sum_{i=1}^{k} (Y_i(t) - \bar{Y}_n(t))\hat{v}_r(t)dt\bigg|\\
&\leq \|\hat{v}_r\|\max_{1 \leq  k < n} w(k/n) \Big \| n^{-1/2}\sum_{i=1}^{k} (Y_i - \bar{Y}_n)\Big\|=\cO_P(g(n))
\end{align*}
for all $1\leq r \leq d$. The rate $g(n)$ follows by using the same arguments as in \cite[Lemma 6.2]{torg2014}, i.e.\ relying on stationarity of the innovations and using the symmetry of the test statistic.\footnote{Here, $g(n)$ and  \autoref{ass:assumptionP1} replace the law of iterated logarithm which is used originally in \citet[Theorem 2.8.1]{horvath1997limit}.}  Hence, as a direct consequence we get
\[
	\max_{1 \leq  k < n} w(k/n) |n^{-1/2}\sum_{i=1}^{k}(\hat{\bY}_i - \bar{\hat{\bY}}_n)|_{\hat{\bSigma}}=\hat{\lambda}_d^{-1/2}\cO_P(g(n)).
\] 
for the projected counterpart. Set $m=\ntheta$ as before and observe that 
\begin{equation}\label{eq:sqrt}
	\max_{1 \leq  k < m-\lfloor n\alpha\rfloor}\left(\frac{k}{n-k}\right)
\begin{cases}
= m/(n-m), &\alpha =0,\\
\leq (m-n\alpha)/(n-m), &0<\alpha<\theta.
\end{cases}
\end{equation}
Now, on one hand, by taking \eqref{eq:decomp_ch} and $w(k/n)k/n=(k/(n-k))^{1/2}$  into account and by considering the square root version of \eqref{eq:sqrt}  we get
\begin{align*}
	&\max_{1 \leq  k < m-\lfloor n\alpha\rfloor} w(k/n) | \sum_{i=1}^{k}(\hat{\bX}_i - \bar{\hat{\bX}}_n)|_{\hat{\bSigma}}\\
&\leq \max_{1 \leq  k < m} w(k/n) | \sum_{i=1}^{k}(\hat{\bY}_i - \bar{\hat{\bY}}_n)|_{\hat{\bSigma}}  +  |\hat{\bDelta}|_{\hat{\bSigma}} \max_{1 \leq  k < m-n\alpha} w(k/n)\frac{k}{n} (n-m)\\
&=\cO_P\left(n^{1/2}g(n)/\hat{\lambda}_d^{1/2}\right) +  (n-m)^{1/2}(m-n\alpha)^{1/2}|\hat{\bDelta}|_{\hat{\bSigma}}
\end{align*}
for $\alpha \in [0,\theta)$. Whereas, on the other hand, we arrive at
\begin{align*}
	&\max_{1 \leq  k < m} w(k/n) |\sum_{i=1}^{k}(\hat{\bX}_i - \bar{\hat{\bX}}_n)|_{\hat{\bSigma}}\\
&\geq \bigg|\max_{1 \leq  k < m} w(k/n) |\sum_{i=1}^{k}(\hat{\bY}_i - \bar{\hat{\bY}}_n)|_{\hat{\bSigma}}   -  |\hat{\bDelta}|_{\hat{\bSigma}} \max_{1 \leq  k < m} w(k/n)\frac{k}{n} (n-m)\bigg|\\
&=\bigg|A_1  -  |\hat{\bDelta}|_{\hat{\bSigma}} (n/w(m/n))\bigg|,
\end{align*}
where $A_1=\cO_P\left(n^{1/2}g(n)/\hat{\lambda}_d^{1/2}\right)$, as before. \autoref{ass:assumptionE1} and \autoref{ass:assumptionB2} ensure that
\[
	\frac{ n^{1/2}g(n)/\hat{\lambda}_d^{1/2}}{|\hat{\bDelta}|_{\hat{\bSigma}}(n/w(m/n))}\leq \frac{c}{\big|\int\!\Delta(t)\hat{v}_r(t)dt\big|} \frac{g(n)}{n^{1/2}}\Bigg[\frac{\hat{\lambda}_1}{\hat{\lambda}_d}\Bigg]^{1/2} =o_P(1)
\]
for some $c>0$ since $|\hat{\bDelta}|^{-1}_{\hat{\bSigma}}\geq \hat{\lambda}_1^{-1/2}|\int\!\Delta(t)\hat{v}_j(t)dt|$ for $j=1,\ldots,d$. Altogether, we obtain 
\[
	\frac{\max_{1 \leq  k < m-\lfloor n\alpha\rfloor} w(k/n) |\sum_{i=1}^{k}(\hat{\bX}_i - \bar{\hat{\bX}}_n)|_{\hat{\bSigma}}}{|\hat{\bDelta}|_{\hat{\bSigma}}(n/w(m/n))} \begin{cases} =1 + o_P(1), &\alpha = 0\\
\leq (1-\alpha/\theta)^{1/2} + o_P(1), &0<\alpha<\theta
\end{cases}
\] 
which completes this proof.
\end{proof}
\begin{proof}[Proof of \autoref{thm:null_lpm}]  Property \eqref{eq:kernelL2} is stated in Theorem 1 of \citet{reeder2012twoB} and \autoref{ass:assumptionP2} follows from ergodicity and stationarity. However, note that \autoref{ass:assumptionP2} is also immediately implied by \citet[Theorem 3.3]{berkes2013rice}. We proceed with the verification of \autoref{ass:assumptionL}. Going carefully through the proof of \citet[Theorem 4.1.3]{horvath1997limit} and taking \citet[Theorem 2.1.4]{schmitz2011} into account - replacing all considerations for univariate time series with multivariate analogues - we see that it suffices to show the conditions \eqref{C1}, \eqref{C2} and \eqref{C3} below (cf. also \citet[Theorem 1.2.1]{Kirch2014}). Hereby, it is crucial that \lpm-approximable time series fulfill \autoref{ass:assumptionM} by definition. For one thing, we need an approximation of the projected innovations $\{\bY_i\}$ (see \eqref{eq:proj_innov}) by centered multivariate Brownian motions $\{\bW_1(n)\}$, $\{\bW_2(n)\}$ with covariance matrix $\bSigma$ (cf. \autoref{rem:covmatrix}). More precise, we need that
\begin{align}
	\Big|\sum_{i=1}^n \bY_i - \bW_1(n)\Big| &=\cO(n^{1/2-\eta}) \qquad \text{a.s.},\label{C1}\tag{C1}\\
	\Big|\sum_{i=1}^n \bY_{-i} - \bW_2(n)\Big| &=\cO(n^{1/2-\eta}) \qquad \text{a.s.}\label{C2}\tag{C2}
\end{align} 
for some $\eta>0$. (We do not impose any restriction on the dependence structure between $\{\bW_1(n)\}$ and $\{\bW_2(n)\}$). For another thing, we need the asymptotic independence as well as the exact asymptotic distributions of
\begin{align*}
	A_n^* &:= a(\log n) A_n - b_d(\log n),\\
	B_n^* &:= a(\log n) B_n - b_d(\log n)
\end{align*}
with
\begin{align*}
	A_n &= \max_{1\leq k \leq n/\log n} k^{-1/2}\bigg|\sum_{i=1}^k\bY_i\bigg|,\\
	B_n &= \max_{n-n/\log n \leq k < n} (n-k)^{-1/2}\bigg|\sum_{i=k+1}^n \bY_i\bigg|,
\end{align*}
i.e.\ that
\begin{align}\label{C3}
\begin{split}
	\lim_{n\rightarrow\infty}P(A^*_n\leq s,  B^*_n\leq t) &= \lim_{n\rightarrow\infty}P(A^*_n\leq s)P(B^*_n\leq t)\\
 &= \exp(-\exp(-s))\exp(-\exp(-t))
\end{split}\tag{C3}
\end{align}
holds true for all $s,t \in\mR$.

These conditions \eqref{C1}, \eqref{C2} and \eqref{C3} are analogous to conditions A.3 (i)-(iii) of \citet[Theorem 1.2.1]{Kirch2014}. Condition \eqref{C1} replaces Assumption A.3 (i), \eqref{C2} corresponds to Assumption A.3 (ii)  of \citet{Kirch2014} and condition \eqref{C3} corresponds to Assumption A.3 (iii) (up to normalizing sequences $a(\log n)$ and $b_d(\log n)$, cf. \citet[Theorem 4.1.3]{horvath1997limit}). Notice that the additional Assumption A.1 of \citet{Kirch2014} is trivially fulfilled in our setting, due to the shape of our test statistic.  
\\
\\
We begin by discussing and verifying $\eqref{C1}$. It is easy to see, that the projected time series $\{\bY_{i}\}_{i\in\mZ}$ remains \lpm-approximable (now in $\mR^d$) with same $\kappa>2$ and same rate $\delta(m)$. In particular $E\|Y_0\|^\kappa<\infty$, $\kappa>2$, holds true and $\bSigma$ is obviously positive definite due to $\lambda_d>0$ (cf. also \autoref{rem:covmatrix}). Hence, $\eqref{C1}$ follows immediately from Theorem A.1 of \citet[cf. Theorem S2.1 in the supplement]{Aue2014} taking \citet[disp. (A.1.16)]{horvath1997limit} into account. Furthermore, a careful examination of the proof of \citet[Theorem A.1]{Aue2014} shows that their arguments do not rely on causality and their strong approximation in Theorem A.1 could be restated for general (noncausal) \lpm-approximable multivariate time series using \autoref{def:lpm} with $H=\mR^d$ and $\delta(m)=m^{-\gamma}$ for some $\gamma>2$. Note that in the proof of Theorem A.1 of \citet{Aue2014}, coupling expressions like, e.g.,
\[
   E[\bY_{0,r}\bY_{j,r}]=E\big[\bY_{0,r}(\bY_{j,r}-\bY_{j,r}^{(j)})\big],
\]
are only valid under causality, since $E[\bY_{0,r}\bY_{j,r}^{(j)}]=0$ is necessary. However, relation $\smash{E[\bY^{(j)}_{0,r}\bY_{j,r}^{(j)}]=0}$ is valid under noncausality and it is known that the above expression can be easily replaced by 
\begin{equation}\label{eq:coupling_trick}
 E[\bY_{0,r}\bY_{j,r}]=E\big[\bY_{j,r}(\bY_{0,r}-\bY_{0,r}^{(j)})\big] +  E\big[\bY^{(j)}_{0,r}(\bY_{j,r}-\bY_{j,r}^{(j)})\big].
\end{equation}
Now, observe that after time inversion $\{\bY_{-i}\}_{i\in\mZ}$ remains \lpm-approximable (in the sense of \autoref{def:lpm}) with the same $\kappa$, the same rate $\delta(m)$ and with the same long run covariance matrix $\bSigma$. Hence, according to previous considerations, \eqref{C2} holds true, as well.

Finally, we verify \eqref{C3}. In the setting of linear processes, \citet[Theorem 4.1.3]{horvath1997limit} have shown asymptotic independence of $A_n$ and $B_n$ by replacing the $\bY_i$'s  in $B_n$ by truncated approximations (cf. \citet[p.308]{horvath1997limit}). We adapt this approach in a straightforward manner by considering $m$-dependent copies $\bY_i^{(m)}$ (cf. \eqref{eq:F_mcopy}) and defining
\[
B_n':=\max_{n-n/\log n \leq k < n} (n-k)^{-1/2}|\sum_{i=k+1}^n \bY_i^{(m_n)}|,
\]
where $m_n:=n-3n/\log n$. The representation of $\bY_i$'s as a shift of i.i.d.\ random variables and the construction of the $\bY^{(m_n)}_i$'s  ensures that $Z_{k,r}:=\bY_{k,r}-\bY_{k,r}^{(m_n)}$ are equally distributed for all $k$. Hence, it holds that
\[
	E |Z_{k,r}|^2= E |Z_{0,r}|^2=\cO(m_n^{-2\gamma})
\]
for all $k\in\mZ$ and $r=1,\ldots,d$. Furthermore,
\begin{align*}
	|B_n-B_n'| &\leq \max_{n-n/\log n \leq k < n} (n-k)^{-1/2}\Big|\sum_{i=k+1}^n \left(\bY_i - \bY_i^{(m_n)}\right)\Big|\\
				 &\leq d \sum_{r=1}^d \left(\max_{n-n/\log n \leq k < n} (n-k)^{-1/2}\big|\sum_{i=k+1}^n Z_{i,r}\big|\right).
\end{align*}
An application of the H\'{a}jeck-R\'{e}nyi type inequality of \citet[Theorem 2]{weng1969} yields that
\begin{align*}
&P\bigg(\max_{n-n/\log n\leq k < n}(n-k)^{-1/2}\big|\sum_{i=k+1}^nZ_{i,r}\big|>(\log n)^{-1/2}\bigg)\\
&\leq \bigg((\log n)^{1/2} \sum_{k=n-n/\log n}^{n-1} (n-k)^{-1/2} \big(E |Z_{k,r}|^2\big)^{1/2} \bigg)^2\\
&= (\log n) E |Z_{0,r}|^2\bigg(\sum_{k=1}^{n/\log n} k^{-1/2}\bigg)^2=\cO\left((\log n) m_n^{-2\nu} \frac{n}{\log n}\right)=\cO(n^{1-2\nu}),
\end{align*} 
which implies
\[
B_n=B_n'+o_P(a(\log n)^{-1}).
\]
Now, observe that $B_n'$ and $A_n$ are independent because the sets $\{\bY_i,  i \leq n/\log n\}$ and $\smash{\{\bY_i^{(m_n)}, i\geq  n-n/\log n, \; m_n=n-3n/\log n\}}$ are obviously independent for sufficiently large $n$. Finally, \eqref{C3} follows from \citet[Lemma 2.2]{lajos1993normal} taking \citet[Lemma 29.5]{Davidson1994} into account.
\end{proof}

\begin{proof}[Proof of \autoref{prop:moments2max}]  The well known results of \citet{moricz1976moment} show that moment inequalities for partial sums yield analogous moment inequalities for maxima of partial sums. Furthermore, in \citet{tomacsa2006hajek} it is shown that inequalities for maxima of partial sums yield inequalities for weighted maxima of partial sums and vice versa. Carefully inspecting the proofs of \citet[Theorem 1]{moricz1976moment} and of \citet[Theorem 2.1]{tomacsa2006hajek} we observe that the same results can be restated in our functional setting with $\kappa>2$, as well. Therefore, \citet[Theorem 1]{moricz1976moment} together with assumption \eqref{eq:econd} and Markov's inequality yield
\begin{equation}\label{eq:mza}
	x^{\kappa} P\left(\max_{1\leq  k \leq n} \big\|\sum_{i=1}^k Y_i\big\| \geq x \right) \leq E\bigg[\max_{1\leq k \leq n}\big\|\sum_{i=1}^k Y_i\big\|^{\kappa}\bigg]\leq c_1 n^{\kappa/2}
\end{equation}
for all $x>0$, $n\in \mN$ and some $c_1>0$. Next, we use $n^{\kappa/2}=\cO(\sum_{k=1}^n k^{\kappa/2-1})$ and apply \citet[Theorem 2.1]{tomacsa2006hajek} to obtain,
\[
	x^{\kappa} P\left(\max_{1\leq  k \leq n}k^{-1/2} \big\|\sum_{i=1}^k Y_i\big\| \geq x \right) \leq c_2 \sum_{k=1}^n k^{-1} \leq c_3 \log n,
\]
for all $x>0$, $n\in \mN$ and some $c_2$, $c_3>0$. The conclusion follows on setting $x=c_4(\log n)^{1/\kappa}$ with a suitable constant $c_4>0$.
\end{proof}

\begin{remark}
In the previous proof we applied \citet[Theorem 1]{moricz1976moment} with $g(F_{b,n})=n^\alpha$, $b=0$ and $\alpha=\kappa/2>1$. In case of $\kappa=2$ (i.e.\ $\alpha=1$) an additional logarithmic term would appear on the right-hand side of \eqref{eq:mza}  (cf. \citet[Theorem 3]{moricz1976moment}).
\end{remark}

We proceed with the proof of \autoref{prop:twosided}. \citet[Proposition 4]{berkes2011split} have shown the corresponding result in the univariate setting and their techniques, slightly modified, are directly applicable to the functional setting, as shown below. Here, we demonstrate that Proposition 4 of \citet{berkes2011split} is extensible to {\it{noncausal}} centered \lpm-approximable functional time series $\{Y_i\}_{i\in \mZ}$. We want to emphasize that another, more sophisticated extension - yet for {\it{causal}} centered \lpm-approximable time series -, has been developed by \citet[cf. Theorems 3.1, 3.3 and Remark 3.2]{berkes2013rice}. 

\begin{proof}[Proof of \autoref{prop:twosided}] We want to point out that - up to the functional setting - the proof presented here is for most parts identical with \citet[Proposition 4]{berkes2011split} and that we stay very close to their exposition. To avoid misunderstandings and for the convenience of the reader, we restate their proof in the functional setting, emphasizing the necessary modifications. In adaption to our situation, the major difficulty stems from the relations (37) - (39) of \citet{berkes2011split} which are not clear in the functional setting for arbitrary $\kappa>2$ and therefore are substituted by \eqref{eq:maindecomprel} below. This is done using a result of \citet{berkes2013rice} which is, however, restricted to $\kappa\in(2,3)$.
\\
\\
Let $S_n=\sum_{i=1}^n Y_i$ denote the partial sums, let $\psi_2>0$, $\psi_\kappa>0$ be arbitrary constants and let
\[
	D_\kappa:=\sum_{m=0}^\infty (E\|X_0-X_0^{(m)}\|^\kappa)^{1/\kappa}<\infty,
\]
where the finiteness holds true in view of $\sum_{m=1}^\infty \delta(m)<\infty$. The cases $\kappa=2$ and $2<\kappa<3$ are treated separately, where the former one can be seen as follows: Via the decomposition
\[
		Y_0(t)Y_j(s) =\Big(Y_0(t)-Y_0^{(j)}(t)\Big)Y_j(s) +  Y^{(j)}_0(t)\Big(Y_j(s)-Y_j^{(j)}(s)\Big) + Y_0^{(j)}(t)Y_j^{(j)}(s),
\] 
using stationarity and $D_2<\infty$ we obtain
\begin{equation}\label{eq:case2}
	E\|S_n\|^2=\sum_{k=1}^n E\int\!Y_k(t)Y_k(t)dt + 2\!\!\sum_{1\leq k< l \leq n}  E\int\! Y_k(t)Y_l(t)dt \leq n C_2
\end{equation}
for some $C_2>D_2/\psi_2$ and all $n\in\mN$, which finishes the proof for $\kappa=2$. For more details cf. \citet{berkes2011split}. For the latter case, i.e.\ $2<\kappa<3$, the idea is to show, that for any $n_0\in \mN$ there is some $C_\kappa(n_0)$ such that:
\begin{equation}\label{eq:prop2n}
	E\|S_n\|^\kappa \leq C_\kappa n^{\kappa/2} \quad \forall n \leq n_0\qquad \Rightarrow \qquad E\|S_n\|^\kappa \leq C_\kappa n^{\kappa/2} \quad \forall n \leq 2n_0.
\end{equation}
Hence, by induction, it is possible to conclude that $C_\kappa$ does not depend on $n_0$ which then completes the proof. 
\\
\\
Now, \eqref{eq:prop2n} can be verified by selecting an arbitrary $n_0\in \mN$ and choosing $C_\kappa (n_0)$ large enough, such that on the one hand
\begin{equation}\label{eq:partialsum_upton0}
	E\|S_n\|^\kappa \leq C_\kappa (n_0) n^{\kappa/2}
\end{equation}
for all $n\leq n_0$ and on the other hand $C^{1/\kappa}_\kappa(n_0) > D_\kappa/\psi_\kappa$. Using Jensen's inequality we have 
\[
	E\|Y_k-Y_k^{(n-k)}\|\leq (E\|Y_k-Y_k^{(n-k)}\|^\kappa)^{1/\kappa},
\]
which, via basic inequalities for norms, yields that
\begin{equation}\label{eq:snestmain}
	E\|S_{2n}\|^\kappa\leq \left(2 D_\kappa + \left(E\bigg[\|Z_n + W_n\|^\kappa\bigg]\right)^{1/\kappa} \right)^\kappa
\end{equation}
where
\[
	Z_n = \sum_{k=1}^n Y_k^{(n-k)},\qquad W_n = \sum_{k=1}^n Y_{n+k}^{(k-1)}
\]
(cf. \citet[disp. (36)]{berkes2011split}). Next, observe that
\[
	E\|Z_n\|^p\leq \left(\bigg(E\|S_n\|^p\bigg)^{1/p} + D_p\right)^p
\]
for $p=2$ or $p=\kappa$, respectively, and that the same holds true if we replace $W_n$ by $Z_n$. Hence, in view of \eqref{eq:partialsum_upton0}, we arrive at
\begin{align}\label{eq:snestprop1}
\begin{split}
	E\|Z_n\|^\kappa &\leq n^{\kappa/2}C_\kappa(1+\psi_\kappa)^\kappa, \\
       E\|W_n\|^\kappa &\leq n^{\kappa/2}C_\kappa(1+\psi_\kappa)^\kappa,\\
\end{split}
\end{align}
for $n\leq n_0$. Furthermore, due to \eqref{eq:case2} it holds also that
\begin{align}\label{eq:snestprop2}
\begin{split}
	E\|Z_n\|^2 &\leq nC_2(1+\psi_2)^2,\\
        E\|W_n\|^2 &\leq nC_2(1+\psi_2)^2,
\end{split}
\end{align}
 for all $n\in\mN$ (cf. \citet{berkes2011split}). Recall that $C^{1/2}_2 > D_2/\psi_2$ and that no restriction on $n$ is needed here due to \eqref{eq:case2}. Observe that $Z_n$ and $W_n$ are mean zero and independent. Therefore, by \citet[Lemma 3.1]{berkes2013rice} we have for $2<\kappa<3$ that
\begin{align}\label{eq:maindecomprel}
\begin{split}
	E\|Z_n + W_n\|^\kappa &\leq E\|Z_n\|^\kappa + E\|W_n\|^\kappa\\
 &\quad + E\|Z_n\|^2 (E\|W_n\|^2)^{\kappa/2-1}\\
 &\quad + E\|W_n\|^2 (E\|Z_n\|^2)^{\kappa/2-1}.
\end{split}
\end{align}
Now, combining \eqref{eq:snestprop1}, \eqref{eq:snestprop2} and \eqref{eq:maindecomprel} yields
\begin{equation}\label{eq:snestsecpart}
	E\|Z_n + W_n\|^\kappa \leq 2(1+\psi_\kappa)^\kappa C_\kappa n^{\kappa/2}+ 2 \left((1+\psi_2)^2 C_2 n\right)^{\kappa/2}
\end{equation}
for all $n\leq n_0$, which is a simple but significant modification of \citet[disp. (37)]{berkes2011split}. Consequently,  from \eqref{eq:snestmain} and \eqref{eq:snestsecpart} we obtain
\begin{align*}
E\|S_{2n}\|^\kappa &\leq \left\{2D_\kappa + \left(E\|Z_n + W_n\|^\kappa\right)^{1/\kappa}\right\}^\kappa\\
&\leq  \left\{2D_\kappa + \left(2(1+\psi_\kappa)^\kappa C_\kappa n^{\kappa/2}+ 2 \left((1+\psi_2)^2 C_2 n\right)^{\kappa/2}\right)^{1/\kappa}\right\}^\kappa\\
&\leq  C_\kappa (n\Psi)^{\kappa/2}
\end{align*}
where
\[
	\Psi(\psi_\kappa, \psi_2,C_\kappa,C_2):=2\psi_\kappa + \left[2(1+\psi_\kappa)^\kappa + 2(1+\psi_2)^\kappa C_2^{\kappa/2}C_\kappa^{-1}\right]^{1/\kappa}\downarrow 2^{1/\kappa}<2
\]
as $(\psi_2,\psi_\kappa,C_\kappa^{-1})\rightarrow 0$ (cf. \citet[disp. (38)]{berkes2011split}). Copying the final arguments of \citet[Proof of Proposition 4]{berkes2011split} completes the proof.
\end{proof}

Next, we take a closer look at the estimation of the eigenstructure of $\cC$. 

\begin{proof}[Proof of \autoref{thmp2:estimation_longrun}] Due to the symmetry of $K(x)$ the estimator $\hat{\zeta}$ can be rewritten as
\begin{equation}\label{eq:basic_decomposition_estimator}
\hat{\zeta}=\hat{\zeta}_0(t,s) + \sum_{i=1}^{n}K(i/h_n)\hat{\zeta}_i(t,s) + \sum_{i=1}^{n}K(i/h_n)\hat{\zeta}_i(s,t).
\end{equation}
Note that in view of \citet[Theorem 3.1]{hoermann2010weakly} the covariance estimation is of order 
\begin{equation}\label{eq:sm001}
	\iint(\hat{\zeta}_0(t,s)-E[Y_0(t)Y_0(s)])^2dtds=\cO_P(n^{-1}).
\end{equation}
(Their arguments carry over to our case of noncausality in a straightforward manner using modifications similar to \eqref{eq:coupling_trick}.) It remains to investigate the long run part of the estimate, where due to symmetry it suffices to consider the second term of \eqref{eq:basic_decomposition_estimator}.
We define a centered version of the second expression of \eqref{eq:basic_decomposition_estimator}
\[
	\hat{c}_{1,n}(t,s)=\sum_{i=1}^{n}K(i/h_n)\hat{\gamma}_i(t,s)
\]
with $\hat{\gamma}_i (t,s)= n^{-1}\sum_{j=1}^{n-i} Y_{j}(t) Y_{j+i}(s)$ and take into account that the difference between the original expression and its centered counterpart is of order 
\begin{equation}\label{eq:sm002}
\iint \left(\sum_{i=1}^{n}K(i/h_n)\hat{\zeta}_i(t,s) - \hat{c}_{1,n}(t,s)\right)^2 dtds= \cO_P(h_n^2/n)
\end{equation}
(cf. \citet[proof of Theorem 2]{reeder2012twoB}, as before, with straightforward modifications in view of noncausality). The centered version can be decomposed as follows:
\[
	\hat{c}_{1,n}-c_1 = \left[\hat{c}_{1,n}-\hat{c}^{(m_n)}_{1,n}\right] +  \left[\hat{c}^{(m_n)}_{1,n}-E \hat{c}^{(m_n)}_{1,n}\right] + \left[E \hat{c}^{(m_n)}_{1,n} - c^{(m_n)}_1\right] +  \left[c^{(m_n)}_1-c_1\right],
\]
where
\begin{align*}
c_1(t,s) &= \sum_{i=1}^{\infty}E[Y_{0}(t)Y_i(s)],\\
c^{(m_n)}_1(t,s) &= \sum_{i=1}^{m_n}E[Y^{(m_n)}_{0}(t) Y^{(m_n)}_i (s)],\\
\hat{c}^{(m_n)}_{1,n}(t,s) &= \sum_{i=1}^{n}K(i/h_n)\hat{\gamma}^{(m_n)}_i(t,s),\\
\hat{\gamma}^{(m_n)}_i(t,s) &= n^{-1}\sum_{j=1}^{n-i} Y^{(m_n)}_{j}(t) Y^{(m_n)}_{j+i}(s)
\end{align*}
and $(m_n)$ indicates the $m_n$-dependent versions. The sequence $m_n$ needs to fulfill $m_n=o(h_n)$ and $m_n\rightarrow \infty$, as $n\rightarrow\infty$. The main extension of the proof of \citet{reeder2012twoB} is the introduction of the additional term $E[\hat{c}^{(m_n)}_{1,n}(t,s)]$ and that we allow for an increase in the dependency of $\hat{c}^{(m_n)}_{1,n}(t,s)$ with increasing sample size $n\rightarrow\infty$. We proceed by observing that
\begin{align}\label{varest}
\begin{split}
	&\bigg|\Var\bigg(\sum_{i=1}^{n}K(i/h_n)\hat{\gamma}^{(m_n)}_i(t,s)\bigg)\bigg|\\
&= \bigg|\sum_{i,j=1}^{n}\Cov\bigg(K(i/h_n)\hat{\gamma}^{(m_n)}_i(t,s),K(j/h_n)\hat{\gamma}^{(m_n)}_j(t,s)\bigg)\bigg|\\
	&\leq c_1\sum_{i,j=1}^{n}\left|\Cov\left(\hat{\gamma}^{(m_n)}_i(t,s),\hat{\gamma}^{(m_n)}_j(t,s)\right)\right|\\
 &\leq n^{-2}c_1\sum_{i,j=1}^{n}\sum_{k,r=1}^n\left|\Cov\left(Y^{(m_n)}_{k}(t) Y^{(m_n)}_{k+i}(s),Y^{(m_n)}_{r}(t) Y^{(m_n)}_{r+j}(s)\right)\right|
\end{split}
\end{align}
for some $c_1>0$. Hence, \eqref{varest}, stationarity and $m_n$-dependence yield, by counting the independent terms and taking into account that $K(x)\equiv0$ for $x>c$ for some $c>0$,
\[
n^{2}\iint \bigg|\Var\left[\sum_{i=1}^{n}K(i/h_n)\hat{\gamma}^{(m_n)}_i(t,s)\right]\bigg|dtds=\cO(1)\sum_{i,j=1}^{\lfloor ch_n \rfloor}\sum_{k,r=1}^n \delta_{k,r}^{i,j},
\]
where 
\[
	\delta_{k,r}^{i,j}:=
\begin{cases} 0,&  r - (k+i) \geq  m_n, r\geq k,\\
0,&  k - (r+j) \geq  m_n, r\leq k,\\
1,&  r - (k+i) \leq  m_n, r\geq k,\\
1,&  k - (r+j) \leq  m_n, r\leq k,
\end{cases}\quad
=~\begin{cases} 1,&  0 \leq r - k \leq m_n + i,\\
1,& 0 \leq k - r \leq m_n + j,\\
0,& \text{else}. 
\end{cases}
\]
Due to stationarity, the values $\delta_{k,r}^{i,j}$ depend only on $i$, $j$ and on the difference of $k-r$. Hence, we have
\begin{align*}
\sum_{i,j=1}^{\lfloor ch_n \rfloor}\sum_{k,r=1}^n \delta_{k,r}^{i,j} &= \sum_{i,j=1}^{\lfloor ch_n \rfloor}\sum_{z=1}^n \left\{\sum_{q=0}^n \delta_{z,z+q}^{i,j} + \sum_{q=0}^n \delta_{z+q,z}^{i,j} \right\}\\
&= \sum_{i,j=1}^{\lfloor ch_n \rfloor}\sum_{z=1}^n \left\{\sum_{q=0}^n \delta_{0,q}^{i,j} + \sum_{q=0}^n \delta_{q,0}^{i,j} \right\}\\
&\leq  n\sum_{i,j=1}^{\lfloor ch_n \rfloor}(m_n+i+j)\leq 3 c' n h_n^3
\end{align*}
for some $c'>0$. From above considerations we obtain 
\begin{align}\label{eq:variance_estimation}
\begin{split}
	&E\iint \left(\hat{c}^{(m_n)}_{1,n}(t,s)-E\hat{c}^{(m_n)}_{1,n}(t,s)\right)^2dtds\\
&=\iint \Var[\hat{c}^{(m_n)}_{1,n}(t,s)]dt ds =\cO(h^3/n).
\end{split}
\end{align}

 Next, using standard arguments and stationarity we see that
\begin{equation}\label{sm1}
  \begin{split}
	& \left(\iint \left\{  E \hat{c}^{(m_n)}_{1,n}(t,s)-c^{(m_n)}_1(t,s)\right\}^2 dtds\right)^{1/2} \\
 &= \left(\iint \left\{ \sum_{i=1}^{m_n}\big(K(i/h_n)(n-i)/n-1\big)E[Y^{(m_n)}_{0}(t) Y^{(m_n)}_{i}(s)]\right\}^2 dtds\right)^{1/2} \\
 &= \cO\bigg(n^{-1}\sum_{i=1}^{m_n}i\left(\iint \left\{ E[Y^{(m_n)}_{0}(t) Y^{(m_n)}_{i}(s)]\right\}^2 dtds\right)^{1/2} \\
&+ \bigg\{\max_{i=1,\ldots,m_n}|K(i/h_n)-1|(h_n/i)^\rho\bigg\} \\
&\times \bigg\{h_n^{-\rho}\sum_{i=1}^{m_n}i^\rho\left(\iint \left\{ E[Y^{(m_n)}_{0}(t) Y^{(m_n)}_{i}(s)]\right\}^2 dtds\right)^{1/2}\bigg\}\bigg)
\\
&=\cO\bigg(h_n^{-\rho}\sum_{i=1}^{m_n}i^\rho\left(\iint \left\{ E[Y_{0}(t) Y_{i}(s)]\right\}^2 dtds\right)^{1/2}+m_n/\exp(c m_n)\bigg)
\end{split}
\end{equation} 
for some $c>0$. The last line follows since $m_n=\cO(h_n)=\cO(n)$ and by decomposing as follows
\begin{align}\label{decompos}
\begin{split}
	&Y^{(m)}_{0}(t) Y^{(m)}_{i}(s)-Y_0(t) Y_{i}(s)\\
&= Y^{(m)}_{0} (t) \left(Y^{(m)}_{i}(s)-Y_{i}(s)\right) + \left(Y^{(m)}_{0}(t)-Y_0(t)\right) Y_{i} (s).
\end{split}
\end{align}
Now, by \citet[proof of Theorem 2]{reeder2012twoB} and the exponential decay of $\delta(m)$ we observe that
\begin{align}\label{sm1a}
  \begin{split}
	E\|\hat{c}_{1,n}-\hat{c}^{(m_n)}_{1,n}\| &= \cO\Big( m_n \left\{E\|Y_0-Y_0^{(m_n)}\|^2\right\}^{1/2}\\ 
&\quad+ \sum_{i=m_n+1}^\infty\left\{E\|Y_0-Y_0^{(i)}\|^2\right\}^{1/2} \Big) = \cO(m_n/\exp(c m_n))
\end{split}
\end{align}
for some $c>0$. Using decomposition \eqref{decompos}, stationarity and again the exponential decay of $\delta(m)$ we get
\begin{align}\label{eq:cov_m_dep}
\begin{split}
	\|c^{(m_n)}_1-c_1\| &\leq  \left(\iint \left(\sum_{i=1}^{m_n}E[Y^{(m_n)}_{0}(t)Y^{(m_n)}_{i}(s)-Y_0(t) Y_{i}(s)]\right)^2 dtds\right)^{1/2}\\
&\quad + \left(\iint \left(\sum_{i= m_n+1}^\infty E[Y_{0}(t) Y_{i}(s)]\right)^2 dtds\right)^{1/2}\\
&= \cO\left(m_n \left\{E\|Y_0-Y_0^{(m_n)}\|^2\right\}^{1/2}\right)\\
&\quad + \bigg(\iint \bigg(\sum_{i= m_n+1}^\infty E\Big[\big(Y_{0}(t)-Y^{(i)}_0(t)\big)Y_i(s)\\
&\quad +Y^{(i)}_{0}(t)\big(Y_{i}(s)-Y^{(i)}_i(s)\big)\Big]\bigg)^2 dtds\bigg)^{1/2} \\
&= \cO(m_n/\exp(cm_n))
\end{split}
\end{align}
for some $c>0$. Combining \eqref{eq:basic_decomposition_estimator} - \eqref{eq:cov_m_dep}, we get  
\[
	\|\hat{\zeta}-\zeta\|=\cO_P\left( (h_n/n)^{1/2} h_n + h_n^{-\rho}\sum_{i=1}^{m_n}i^\rho \delta(i)+m_n/\exp(c m_n)\right).
\]
Setting $m_n:=\lfloor (\log n) / c \rfloor$ the last term becomes negligible (in comparison to the first term) and we obtain the desired rate. 
\end{proof}

{\footnotesize   
\begingroup
\setstretch{0.1}
\microtypesetup{tracking=false}

\endgroup
}


\begin{thebibliography}{50} 

		\bibitem[{Aston \& Kirch(2012)}]{kirch2012functional}
		{Aston J.A.D., Kirch C.} (2012)  \href{http://dx.doi.org/10.1016/j.jmva.2012.03.006}{Detecting and estimating changes in dependent functional data. \journalmultvaranal, 109:204--220 }  


		\bibitem[{Aue \textit{et~al.\!}(2009)H{\"o}rmann, Horv{\'a}th and Reimherr}]{Aue2009}
		{Aue A., H{\"o}rmann S., Horv\!\!\;{\'a}th L., Reimherr M.} (2009)  \href{http://dx.doi.org/10.1214/09-AOS707}{Break detection in the covariance structure of multivariate time series models. \annstat, 37(6B):4046--4087 }

\bibitem[{Aue \textit{et~al.\!}(2014)H{\"o}rmann, Horv{\'a}th and Hu{\v{s}}kov{\'a}}]{Aue2014} 
		{Aue A., H{\"o}rmann S., Horv\!\!\;{\'a}th L., Hu{\v{s}}kov\!\!\;{\'a} M.} (2014) \href{http://dx.doi.org/10.5705/ss.2012.233}{Dependent functional linear models with applications to monitoring structural change. \statisticasinica{: Supplement}, 24(3):S1--S13}

		\bibitem[{Aue \& Horv{\'a}th(2013)}]{aue2013}
		{Aue A.,  Horv\!\!\;{\'a}th L.} (2013) \href{http://dx.doi.org/10.1111/j.1467-9892.2012.00819.x}{Structural breaks in time series. \journaltimeseranal, 34(1):1--16}

\bibitem[{Aue \textit{et~al.\!}(2015)}]{aue2015}
    	Aue A., Rice G., S{\"o}nmez O. (2015) \href{http://arxiv.org/abs/1511.04020v1}{Dating structural breaks in functional data without dimension reduction. \arxiv: 1511.04020v1, 1--30. Preprint.}

\bibitem[{Berkes \textit{et~al.\!}(2009) Berkes, Gabrys, Horv{\'a}th, and Kokoszka (2009)}]{berkes2009functional}
		{Berkes I., Gabrys R., Horv\!\!\;{\'a}th L., Kokoszka P.} (2009)   \href{http://dx.doi.org/10.1111/j.1467-9868.2009.00713.x}{Detecting changes in the mean of functional observations. \journalroyalsociety, 71(5):927--946}

		\bibitem[{Berkes \textit{et~al.\!}(2013)Berkes, Horv{\'a}th, and Rice}]{berkes2013rice}
		{Berkes I., Horv\!\!\;{\'a}th L., Rice G.} (2013) \href{http://dx.doi.org/10.1016/j.spa.2012.10.003}{Weak invariance principles for sums of dependent random functions. \stochprocappl, 123(2):385--403}
		 
     	\bibitem[{Berkes \textit{et~al.\!}(2015)Berkes, Horv{\'a}th, and Rice}]{berkes2015rice}
		Berkes I., Horv\!\!\;{\'a}th L., Rice G. (2015)   \href{http://arxiv.org/abs/1503.00741v2}{On the asymptotic normality of kernel estimators of the long run covariance of functional time series. \arxiv: 1503.00741v2, 1--37. Preprint, second version.}

		\bibitem[{Berkes \textit{et~al.\!}(2011)Berkes, H{\"o}rmann, and Schauer}]{berkes2011split}
		{Berkes I., H{\"o}rmann S., Schauer J.} (2011)  \href{http://dx.doi.org/10.1214/10-AOP603}{Split invariance principles for stationary processes. \annprob, 39(6):2441--2473}


		\bibitem[{Cs{\"o}rg{\H o} \& Horv{\'a}th(1997)}]{horvath1997limit}
		{Cs{\"o}rg{\H o} M., Horv\!\!\;{\'a}th L.} (1997) {Limit Theorems in Change-Point Analysis}. Wiley, Chichester.

\bibitem[{Chochola \textit{et~al.\!}(2013)Chochola, Hu{\v{s}}kov{\'a}, Pr{\'a}{\v{s}}kov{\'a} and Steinebach}]{chochola2013}
		{Chochola O., Hu{\v{s}}kov\!\!\;{\'a} M., Pr{\'a}{\v{s}}kov\!\!\;{\'a} Z., Steinebach, J.G.} (2013)  \href{http://dx.doi.org/10.1016/j.jmva.2012.10.019}{Robust monitoring of CAPM portfolio betas. \journalmultvaranal, 115:374--395}

\bibitem[{Davidson(1994)}]{Davidson1994}
Davidson J. (1994) Stochastic Limit Theory: An Introduction for Econometricians. Oxford university press, New York.

	\bibitem[{Gabrys \& Kokoszka(2007)}]{gabrys2007}
		{Gabrys R., Kokoszka P.} (2007)  \href{http://dx.doi.org/10.1198/016214507000001111}{Portmanteau test of independence for functional observations. \jasa, 102(480):1338--1348}

\bibitem[{Gombay \& Horv{\'a}th(1996)}]{gombay1996rate}
		{Gombay E., Horv\!\!\;{\'a}th L.} (1996)  \href{http://dx.doi.org/10.1006/jmva.1996.0007}{On the rate of approximations for maximum likelihood tests in change-point models. \journalmultvaranal, 56(1):120--152}



\bibitem[{H{\"o}rmann \& Kidzi{\'n}ski(2015)}]{1lukasz2014} 
{H{\"o}rmann S., Kidzi{\'n}ski \L.} (2015)  \href{http://dx.doi.org/10.1111/sjos.12094}{A note on estimation in Hilbertian linear models.  \scandinavian, 42(1):43--62}

\bibitem[{H{\"o}rmann \& Kokoszka(2010)}]{hoermann2010weakly}
{H{\"o}rmann S., Kokoszka P.} (2010) \href{http://dx.doi.org/10.1214/09-AOS768}{Weakly dependent functional data. \annstat, 38(3):1845--1884}

\bibitem[Horv{\'a}th(1993)]{lajos1993normal}
{Horv\!\!\;{\'a}th L.} (1993) \href{http://dx.doi.org/10.1214/aos/1176349143}{The maximum likelihood method for testing changes in the parameters of normal observations. \annstat, 21(2):671--680}

\bibitem[{Horv{\'a}th \& Kokoszka(2012)}]{horvath2012fda}
{Horv\!\!\;{\'a}th L., Kokoszka P.} (2012)  \href{http://dx.doi.org/10.1007/978-1-4614-3655-3}{{Inference for Functional Data with Applications}. Springer Series in Statistics. Springer, New York. }

\bibitem[{Horv{\'a}th \textit{et~al.\!}(2011)Horv{\'a}th, Kokoszka, and Reeder}]{reeder2012twoA}  
{Horv\!\!\;{\'a}th L., Kokoszka P., Reeder R.} (2011)  \href{http://arxiv.org/abs/1105.0019v1}{Estimation of the mean of functional time series and a two-sample problem. \arxiv: 1105.0019v1, 1--32. Preprint of \citet{reeder2012twoB}, first version.}

\bibitem[{Horv{\'a}th \textit{et~al.\!}(2013)Horv{\'a}th, Kokoszka, and Reeder}]{reeder2012twoB} 
{Horv\!\!\;{\'a}th L., Kokoszka P., Reeder R.} (2013) \href{http://dx.doi.org/10.1111/j.1467-9868.2012.01032.x}{Estimation of the mean of functional time series and a two-sample problem. \journalroyalsociety, 75(1):103--122}

\bibitem[{Horv{\'a}th \textit{et~al.\!}(2014)Horv{\'a}th, Kokoszka, and Rice}]{horvath2013testing}  
{Horv\!\!\;{\'a}th L., Kokoszka P., Rice G.} (2014)  \href{http://dx.doi.org/10.1016/j.jeconom.2013.11.002}{Testing stationarity of functional time series. \journalofeconometrics, 179(1):66--82}

\bibitem[{Horv{\'a}th \textit{et~al.\!}(1999)Horv{\'a}th, Kokoszka, and  Steinebach}]{steinebach1999mdependent}
{Horv\!\!\;{\'a}th L., Kokoszka P., Steinebach J.G.} (1999) \href{http://dx.doi.org/10.1006/jmva.1998.1780}{Testing for changes in multivariate dependent observations with an application to temperature changes. \journalmultvaranal, 68(1):96--119}
 

\bibitem[{Horv{\'a}th \& Rice(2014)}]{test2014}
{Horv\!\!\;{\'a}th L., Rice G.} (2014)  \href{http://dx.doi.org/10.1007/s11749-014-0368-4}{Extensions of some classical methods in change point analysis. \journaltest, 23(2):219--255}

\bibitem[{Horv{\'a}th \& Rice(2015)}]{rice2014} 
{Horv\!\!\;{\'a}th L., Rice G.} (2015)  \href{http://dx.doi.org/10.1111/jtsa.12095}{Testing equality of means when the observations are from functional time series. \journaltimeseranal, 36(1):84--108}

\bibitem[{Horv{\'a}th \textit{et~al.\!}(2014)Horv{\'a}th, Rice and Whipple}]{whipple2014}
Horv\!\!\;{\'a}th L., Rice G., Whipple S. (2014)  \href{http://dx.doi.org/10.1016/j.csda.2014.06.008}{\\ Adaptive bandwidth selection in the long run covariance estimator of functional time series. \computationalstatisticsanddataanalysis, 1--18} 

\bibitem[{Jirak(2012)}]{jirak2012}
{Jirak M.} (2012)  \href{http://dx.doi.org/10.1016/j.jmva.2012.05.007}{Change-point analysis in increasing dimension. \journalmultvaranal, 111:136--159}

\bibitem[{Jirak(2013)}]{jirak2013}
{Jirak M.} (2013) \href{http://dx.doi.org/10.1016/j.spl.2013.06.014}{On weak invariance principles for sums of dependent random functionals. \statprob, 83(10):2291--2296}


\bibitem[{Kamgaing \& Kirch(2016)}]{Kirch2014}
{Kamgaing J.T., Kirch C.} (2016) Detection of change points in discrete valued time series. {In: Handbook of Discrete Valued Time series. Handbooks of Modern Statistical Methods.} Chapman \& Hall/CRC. 

\bibitem[{Kokoszka(2012)}]{kokoszka2012}
{Kokoszka P.} (2012)  \href{http://dx.doi.org/10.5402/2012/958254}{Dependent functional data. \isrnjournal, 2012:1--30}

\bibitem[{Kounias \& Weng(1969)}]{weng1969}
{Kounias E. G., Weng T.} (1969) \href{http://dx.doi.org/10.1214/aoms/1177697615}{An inequality and almost sure convergence. \annmathstat, 40(3):1091--1093}

\bibitem[{M{\'o}ricz(1976)}]{moricz1976moment}
{M{\'o}ricz F.} (1976) \href{http://dx.doi.org/10.1007/BF00532956}{Moment inequalities and the strong laws of large numbers.  \probtheoryrelatedfields, 35(4):299--314}

\bibitem[{Ramsay \& Silverman(2005)}]{ramsay2005functional}
{Ramsay J., Silverman B.} (2005)  \href{http://dx.doi.org/10.1007/b98888}{{Functional Data Analysis}. Springer, New York.} 

\bibitem[Schmitz(2011)]{schmitz2011}
{Schmitz A.} (2011) \href{http://kups.ub.uni-koeln.de/4224/}{Limit theorems in change-point analysis for dependent data. PhD Thesis, University of Cologne.}

\bibitem[Sharipov \textit{et~al.\!}(2015)]{tewes2015}
Sharipov O., Tewes J., Wendler M. (2015) \href{http://arxiv.org/abs/1412.0446v2}{Sequential block bootstrap in a Hilbert space with application to change point analysis. \arxiv: 1412.0446v2. Preprint.}

\bibitem[{T{\'o}m{\'a}cs and L{\'i}bor(2006)}]{tomacsa2006hajek}
T{\'o}m{\'a}cs T., L{\'i}bor Z. (2006) A {H}{\'a}jek--{R}{\'e}nyi type inequality and its applications. \annmathinform, 33:141--149

\bibitem[{Torgovitski(2014)}]{torg2014} 
{Torgovitski L.} (2014a) \href{http://dx.doi.org/10.1007/s00184-014-0487-7}{A Darling-Erd\H{o}s-type CUSUM-procedure for functional data. \metrika, 78(1):1--27,  Online version published in 2014. Printed version appeared in 2015.}

		\bibitem[{Torgovitski(2014b)}]{torg2014approx} 
		Torgovitski L. (2014b)  \href{http://arxiv.org/abs/1407.3625}{A Darling-Erd\H{o}s-type CUSUM-procedure for functional data II. \arxiv: 1407.3625v1. Preprint.}


\bibitem[{Vostrikova(1981)}]{vostrikova1981disorder}
{Vostrikova L.} (1981)  \href{http://dx.doi.org/10.1137/1126034}{Detection of a {``disorder''} in a {Wiener} process. \theorprobabappl, 26(2):356--362}

\bibitem[{Zhou(2011)}]{zhou2011maximum}
		{Zhou J.} (2011) \href{http://dx.doi.org/10.1016/j.csda.2011.01.003}{Maximum likelihood ratio test for the stability of sequence of Gaussian random processes. \computationalstatisticsanddataanalysis, 55(6):2114--2127}

\end{thebibliography}
\end{document}